\documentclass[a4paper,reqno]{amsart}
\usepackage[english]{babel}
\usepackage[T1]{fontenc} 
\usepackage[utf8]{inputenc}
\allowdisplaybreaks  % permite cortes de página en fórmulas
\usepackage{amssymb}

\usepackage{faktor}
\usepackage[svgnames,x11names,dvipsnames]{xcolor}
\usepackage{tikz-cd}
\tikzset{isom/.style={draw=none,every to/.append style={edge node={node [sloped, allow upside down, auto=false]{$\simeq$}}}}}
\usepackage{bbm}
\usepackage{mathtools}
\usepackage{microtype}
\usepackage{comment}
\usepackage{orcidlink}
\usepackage{hyperref}
\newif\ifdraft
\draftfalse
\usepackage{xstring}
\usepackage[abspath]{currfile}
\usepackage[yyyymmdd]{datetime}

\usepackage[inline]{enumitem}
\setlist[itemize,1]{leftmargin=*}
\setlist[enumerate,1]{leftmargin=*,label=\upshape\bfseries\arabic*.}
\setlist[enumerate,2]{leftmargin=*,align=left,label=\upshape\bfseries(\roman*),widest=iii}
\newlist{primenumerate}{enumerate}{1}
\setlist[primenumerate,1]{leftmargin=*,label=\upshape\bfseries\arabic*$'$.}
\newlist{enumerateprime}{enumerate}{1}
\setlist[enumerateprime,1]{leftmargin=*,label=\upshape\bfseries\arabic*.\phantom{$'$}}

\ifdraft

\makeatletter
\newcommand\oddfoot[1]{\gdef\@oddfoot{\reset@font#1}}
\newcommand\evenfoot[1]{\gdef\@evenfoot{\reset@font#1}}
\makeatother
\oddfoot{%
\parbox{\textwidth}{
\rule{\textwidth}{0.4pt}\par
\tiny\sffamily%
Archivo fuente: \textcolor{Red3}{\guillemotleft\StrSubstitute{\currfilename}{"}{}\guillemotright}.
Última modificación: \textcolor{Red3}{\today:\currenttime}.
}
}
\evenfoot{%
\parbox{\textwidth}{
\rule{\textwidth}{0.4pt}\par
\tiny\sffamily%
Archivo fuente: \textcolor{Red3}{\guillemotleft\StrSubstitute{\currfilename}{"}{}\guillemotright}.
Última modificación: \textcolor{Red3}{\today:\currenttime}.
}
}
\newtheoremstyle{comentario}%
{}%
{}%
{\sffamily}%
{0pt}%
{\upshape}%
{}%
{0pt}%
{\color{Red2}{\sffamily\thmname{#1}\framebox{\small\textbf{\thmnumber{#2}}}}\mdseries\upshape\thmnote{ (#3)} }%
\theoremstyle{comentario}

\theoremstyle{plain}
\RequirePackage{soulutf8}
\sethlcolor{MistyRose}
\usepackage{xparse}
\NewDocumentCommand{\badref}{v}{\textcolor{Firebrick1}{\framebox{\textbf{??}}\hl{\texttt{#1}}}}
\definecolor{labelkey}{rgb}{0.8,0.8,0.8}
\RequirePackage[outer,right,marginal]{showlabels}

\RequirePackage{xstring}
\usepackage[hyperpageref]{backref} % For drafts
\renewcommand*{\backrefalt}[4]{%
\ifcase#1%
\textcolor{Red1}{(No citations)}%
\or%
\textcolor{Yellow3}{(One citation on page #2.)}%
\else%
\textcolor{Green3}{(#1 citations in pages #2.)}%
\fi%
}

\fi

\theoremstyle{theorem}
\newtheorem{theorem}[equation]{Theorem}
\newtheorem{corollary}[equation]{Corollary}
\newtheorem{lemma}[equation]{Lemma}
\newtheorem*{lemma*}{Lemma}
\newtheorem{proposition}[equation]{Proposition}

\theoremstyle{example}
\newtheorem{definition}[equation]{Definition}
\newtheorem{remark}[equation]{Remark}
\newtheorem{example}[equation]{Example}

\numberwithin{equation}{section} %\numberwithin{defn}{section}

\newcommand\id{\operatorname{id}}

\DeclareMathOperator{\chr}{char}
\DeclareMathOperator{\Hom}{Hom}
\DeclareMathOperator{\End}{End}
\DeclareMathOperator{\Aut}{Aut}
\DeclareMathOperator{\Pic}{Pic}
\DeclareMathOperator{\Har}{\mathbb{H}}

\DeclareMathOperator{\Gal}{Gal}
\DeclareMathOperator{\Ann}{Ann}

\DeclareMathOperator{\Spec}{Spec}
\DeclareMathOperator{\Ram}{Ram}

\renewcommand\emptyset{\varnothing}
\DeclarePairedDelimiterXPP\fpws[1]{}{\lbrack\mkern-1.2mu\lbrack}{\rbrack\mkern-1.2mu\rbrack}{}{\ifblank{#1}{\cdot}{#1}}
\DeclarePairedDelimiterXPP\laurent[1]{}{\lparen\mkern-3.5mu\lparen}{\rparen\mkern-3.5mu\rparen}{}{\ifblank{#1}{\cdot}{#1}}
\DeclarePairedDelimiterXPP\pairing[1]{}{\langle}{\rangle}{}{\ifblank{#1}{\cdot}{#1}}
\DeclarePairedDelimiterXPP\gen[1]{}{\langle}{\rangle}{}{\ifblank{#1}{\cdot}{#1}}
\DeclarePairedDelimiterX\braces[1]{\lbrace}{\rbrace}{\ifblank{#1}{\:\cdot\:}{#1}}
     
\newcommand{\A}{\mathbb{A}}

\newcommand{\bb}{\mathcal{B}}
\newcommand{\C}{\mathbb{C}}
\newcommand{\cc}{\mathcal{C}}
\newcommand{\f}{\mathbbm{f}}
\newcommand{\G}{\mathbb{G}}
\renewcommand{\H}{\mathbb{H}}

\newcommand{\I}{\mathbb{I}}
\renewcommand{\k}{\Bbbk}

\newcommand{\m}{\mathfrak{m}}
\newcommand{\oo}{\mathcal{O}}
\renewcommand{\t}{\mathbbm{t}}
\renewcommand{\P}{\mathbb{P}}
\renewcommand{\ss}{\mathcal{S}}
\newcommand{\Z}{\mathbb{Z}}
\let\stdphi\phi
\let\phi\varphi
\let\varphi\stdphi

\newcommand{\eqdef}{\coloneqq}
\newcommand{\defeq}{\eqqcolon}
\newcommand{\suchthat}{:}
\newcommand{\gkummer}[3][n]{#2\lbrace\mkern0mu{#3}^{1/#1}\mkern0.5mu\rbrace}

\newcommand{\Axn}[2][n]{A_{x}\lbrace\mkern0mu{#2}^{1/#1}\mkern0.5mu\rbrace}
\newcommand{\Kxn}[2][n]{K_{x}\lbrace\mkern0mu{#2}^{1/#1}\mkern0.5mu\rbrace}
\newcommand{\AXn}[2][n]{\A_{X}\lbrace\mkern0mu{#2}^{1/#1}\mkern0.5mu\rbrace}

\newcommand{\AXSn}[3][n]{\A_{X,#3}\lbrace\mkern0mu{#2}^{1/#1}\mkern0.5mu\rbrace}

\newcommand{\Zp}{\mathbb{Z}/(p)}
\newcommand{\Zps}{(\mathbb{Z}/(p))^{*}}
\newcommand{\Zpsf}{\left(\faktor{\Z}{(p)}\right)^{\!*}}

\newcommand{\joinrelshort}{\mathrel{\mkern-8mu}}
\newcommand{\shortlongrightarrow}{\relbar\joinrelshort\rightarrow}
\newcommand{\isomto}{\mathrel{\mathop{\setbox0\hbox{$\mathsurround0pt\shortlongrightarrow$}\ht0=0.3\ht0\box0}\limits^{\hspace{-1pt}\scalebox{1.0}{$\sim$}\mkern2mu}}}

\newcommand{\correspondence}{\overset{1:1}\leftrightsquigarrow}

\newcommand*{\prodprime}{\operatornamewithlimits{%
\mathchoice
{\prod\nolimits\raisebox{1.618ex}{\hspace{0em}\makebox[0pt]{$\scriptstyle\prime$}}\hspace{-0.2em}}%display
{\prod\nolimits\raisebox{1ex}{\hspace{-0.15em}\makebox{$\scriptstyle\prime$}}\hspace{0.1em}}%text
{\prod\nolimits\raisebox{0.618ex}{\hspace{-0.16em}\makebox{$\scriptstyle\prime$}}\hspace{0.1em}}%script
{\prod\nolimits\raisebox{0.3ex}{\hspace{-0.2em}\makebox{$\scriptstyle\prime$}}\hspace{0em}}%scriptscript
}}

\title[Kummer theory over the geometric adeles]{%
Kummer theory over the geometric adeles of an algebraic curve
}%
\author[L.~M.~Navas~Vicente]{Luis Manuel Navas Vicente \orcidlink{0000-0002-5742-8679}}
\email{navas@usal.es}

\author[F.~J.~PLaza~Mart\'in]{Francisco J. Plaza Mart\'in}
\email{fplaza@usal.es}

\address{Departamento de Matem\'aticas and IUFFyM, Universidad de
Salamanca,  Plaza de la Merced 1-4
        \\
        37008 Salamanca. Spain.
		\\
}

\author[Á. Serrano~Holgado]{Álvaro Serrano Holgado}
\email{Alvaro\_Serrano@usal.es}

\thanks{Research of the first author supported by grant PID2021-124332NB-C22 from the Ministerio de Ciencia e Innovación (Government of Spain). \begin{enumerate*}[label=Author \arabic* ORCID:] \item 0000-0002-5742-8679, \item 0000-0002-5532-7567, \item 0000-0003-1695-9147. \end{enumerate*}
}
\address{Departamento de Matem\'aticas, Universidad de
Salamanca,  Plaza de la Merced 1-4
        \\
        37008 Salamanca. Spain.
}

\begin{document}
%%%%%%%%%%%%%%%%%%%%%%%%%%%%%%%%%%%%%%%%%%%%%%%%%%%%%%%%%%%%%%%%%%%%%%%%%%%%%%%%%%%
% OPCIONES LATEX GENERALES
%%%%%%%%%%%%%%%%%%%%%%%%%%%%%%%%%%%%%%%%%%%%%%%%%%%%%%%%%%%%%%%%%%%%%%%%%%%%%%%%%%%
\emergencystretch 3em 
%%%%%%%%%%%%%%%%%%%%%%%%%%%%%%%%%%%%%%%%%%%%%%%%%%%%%%%%%%%%%%%%%%%%%%%%%%%%%%%%%%%
%%%%%%%%%%%%%%%%%%%%%%%%%%%%%%%%%%%%%%%%%%%%%%%%%%%%%%%%%%%%%%%%%%%%%%%%%%%%%%%%%%%

\begin{abstract}
	
Our goal is to give a purely algebraic characterization of finite abelian Galois covers of a complete, irreducible, non-singular curve $X$ over an algebraically closed field $\k$. To achieve this, we make use of the Galois theory of commutative rings, in particular the Kummer theory of the ring of geometric adeles $\A_{X}$.

After we establish the triviality of the Picard group $\Pic(\A_{X})$, the general Kummer sequence for Kummerian rings leads to a characterization of $p$-cyclic extensions of $\A_{X}$ in terms of the closed points of $X$. This is an example of a general local-global principle which we use throughout, allowing us to avoid needing the full spectrum of $\A_{X}$. We prove the existence of primitive elements in $p$-cyclic extensions of $\A_{X}$, which yields explicit invariants lying in $\bigoplus_{x \in X} \Zp$ (summing over closed points) classifying them.

From a group-theoretical point of view, we give a complete characterization of which $p$-cyclic subgroups of the full automorphism group of a given $p$-cyclic extension of $\A_{X}$ endow it with a Galois structure. The result is a stratification by the algebraic ramification of the extension modulo a notion of conjugation or twisting of Galois structures, yielding other invariants, in the form of finite tuples over ramified points, which are related to the previous ones in terms of the local Kummer symbols.

With these results in hand, a forthcoming paper will identify, inside the set of $p$-cyclic extensions of $\A_{X}$, those arising from extensions of the function field of the curve $X$, eventually leading to the algebraic characterization of abelian covers of $X$.

\ifdraft
\leavevmode
\begin{description}

\item[\sffamily Archivo fuente] {\sffamily\textcolor{Red3}{\currfilename}}

\item[\sffamily Última modificación] {\sffamily\textcolor{Red3}{\today:\currenttime}.}

\end{description}
\fi
\end{abstract}

\maketitle

% \footnotesize
% \tableofcontents
% \normalsize

%%%%%%%%%%%%%%%%%%%%%%%%%%%%%%%%%%%%%%%%%%%%%%%%%%%%%%%%%%%%%%%%%%%%%%%%%%%%%%%%%%%
%%%%%%%%%%%%%%%%%%%%%%%%%%%%%%%%%%%%%%%%%%%%%%%%%%%%%%%%%%%%%%%%%%%%%%%%%%%%%%%%%%%
%%%%%%%%%%%%%%%%%%%%%%%%%%%%%%%%%%%%%%%%%%%%%%%%%%%%%%%%%%%%%%%%%%%%%%%%%%%%%%%%%%%
\section{Introduction}
\label{sec:introduction}
% (fold)
%%%%%%%%%%%%%%%%%%%%%%%%%%%%%%%%%%%%%%%%%%%%%%%%%%%%%%%%%%%%%%%%%%%%%%%%%%%%%%%%%%%
%%%%%%%%%%%%%%%%%%%%%%%%%%%%%%%%%%%%%%%%%%%%%%%%%%%%%%%%%%%%%%%%%%%%%%%%%%%%%%%%%%%
%%%%%%%%%%%%%%%%%%%%%%%%%%%%%%%%%%%%%%%%%%%%%%%%%%%%%%%%%%%%%%%%%%%%%%%%%%%%%%%%%%%

This paper is the first in a planned series whose aim is the algebraic characterization of finite abelian Galois covers of a complete, irreducible, non-singular curve $X$ over an algebraically closed field $\k$. By ``algebraic'' we mean that from the very beginning we can avoid the analytic and topological machinery that is part of the standard approach, e.g., the Riemann Existence Theorem, the theory of Riemann surfaces and the fundamental group.

To achieve this, we use the Galois theory of commutative rings, introduced by Auslander and Goldman in~\cite{AuslanderGoldman} and developed by Chase, Harrison and Rosenberg in~\cite{ChaseHarrisonRosenberg} to study finite abelian Galois extensions of the geometric adele ring $\A_{X}$ of the curve. To avoid excessive technicalities, we limit ourselves here to the prime cyclic case in the spirit of~\cite{Borevich}. This involves, as the title indicates, determining the corresponding Kummer theory of $\A_{X}$, which we accomplish employing the updated methods outlined by Greither in~\cite{Greither}.

The second step in our strategy is the characterization, inside the set of such extensions of $\A_{X}$, of those coming from finite abelian Galois extensions of the function field $\Sigma$ of the curve via a natural correspondence, thus leading to a characterization of covers of $X$. We give a brief outline of this future work in \S\ref{sec:future work}. A forthcoming paper will present the complete solution to this problem in the cyclic case.

Recall that the adele ring of $X$ is the restricted direct product
\begin{equation*}
\label{E:adele intro}
	      \mathbb{A}_{X}:= \prodprime_{x \in X} (K_{x},A_{x})
	= \Bigl\{ (\alpha_{x})_{x\in X} \in \prod_{x} K_{x} \suchthat \alpha_{x}\in A_{x}  \text{ for almost all } x \Bigr\}
\end{equation*}
where the product is over the closed points of $X$, $A_{x}$ denotes the completion of the local ring at $x \in X$, and $K_{x}$ is its field of fractions. 

The Kummer theory of rings leads us to focus our attention on $\A_{X}$-algebras of the form
\begin{equation}
\label{E:AXt-intro}
	\AXn{\t} \eqdef \A_{X}[T] / (T^{n} - \t)
\end{equation}
which in the context of this paper we will simply call ``adelic algebras'', and whose general definition is given in \S\ref{subsec:definition of adelic algebra}. We will be interested in the case where $\t$ is an idele. Throughout, $n$ is taken coprime to $\chr(\k)$, and later on $n$ will be assumed for convenience to be a prime $p$.

In general, the properties of a ring extension are related to those of its localizations. For example, characterizing when an extension is Galois (\cite[Theorem 12.2.9. (6)]{Ford} and~\cite{Paques}) may be done in this way. However, because of the rather intricate structure of the maximal spectrum of $\A_{X}$ (see~\cite{SerranoHolgado}), this approach unduly complicates matters and will be seen to be unnecessary, as the particular structure of $\A_{X}$ means we can restrict to the closed points $x \in X$. We shall refer to this feature as a ``local-global'' principle, namely, criteria which relates a global property of an object over $\A_{X}$ to the corresponding property at each closed point $x$ of the curve $X$.

Let us give a brief summary of the paper. Before taking on the Kummer theory of $\A_{X}$ we need a few prior technical results, such as the separability of $\AXn{\t}$ (Proposition~\ref{P:Axt-separable}) and the characterization of isomorphisms between adelic algebras of this form (Theorem~\ref{T:AXt isomorphic iff ram equal and idele ratio}).

Turning to the Galois theory of rings, in Theorem~\ref{T:galois local global} we establish a general local-global principle as described above, characterizing $p$-cyclic extensions of $\A_{X}$ in terms of the closed points of $X$. This relies on the triviality of the Picard group $\Pic(\A_{X})$ (Theorem~\ref{T:Pic adeles}).

The characterization of such extensions is reminiscent of classical Kummer theory for fields, insofar as we prove the existence of primitive elements on which the cyclic group  $\cc_{p}$ acts via a suitably defined character $\chi$ (Theorem~\ref{T:primitive element theorem p-cyclic}), but some care needs to be taken  since $\A_{X}$ is a large ring with many idempotents and zero-divisors.
%Fortunately we can once again develop several local-global criteria which allow us to study such properties over $K_{x}$ or $A_{x}$ at each closed point $x$.

The existence of primitive elements in commutative ring extensions has been in itself the subject of previous work (e.g. \cite{Nagahara,Paques}). Our characterization of $p$-cyclic extensions of $\A_{X}$ (Theorem~\ref{T:classification Cp-extensions AX general}) is made explicit via their use, in the sense that isomorphism classes of extensions correspond to invariants lying in $\bigoplus_{x} \Zp$, the sum being over closed points. As can be seen in Example~\ref{EX:ramified cover P1}, our invariants have a clear geometric significance, in accordance with our overall intention of classifying covers of $X$.

One of the features of the Galois theory of rings, in contrast with the special case of fields, is that a fixed extension can be Galois with respect to many different finite subgroups of its full automorphism group $\G$, which can be a quite large infinite group. To our knowledge, the problem of characterizing such subgroups has not been considered before in the literature. It involves a notion of conjugation of Galois extensions with different but isomorphic Galois groups (Definition~\ref{D:conjugate extensions}). When the groups are in fact the same, this corresponds to ``twisting'' the Galois structure. Galois $p$-cyclic extensions of $\A_{X}$ up to conjugacy can again be explicitly classified by invariants (Theorem~\ref{T:correspondence conjugacy p-cyclic}).

For the adele ring $\A_{X}$, we are able to solve the problem completely in the $p$-cyclic case, characterizing those $p$-cyclic subgroups $G$ of $\G$ which are Galois groups of a given extension (Proposition~\ref{P:pointwise transitive equivalences}). 

More importantly, a completely group-theoretical characterization of conjugacy is given in Theorem~\ref{T:conjugacy subgroups equivalences}. This can be regarded as a version of the classical Galois correspondence, taking into account the above non-uniqueness phenomena.

This characterization involves the algebraic ramification of the extension, which will be linked to the geometric concept in the upcoming work mentioned above showing how to obtain covers of $X$ from Galois extensions of $\A_{X}$. In the present context, we obtain a stratification by ramification (Corollary~\ref{C:stratification conjugacy p-cyclic ramification}) and another set of invariants, this time consisting of finite tuples over ramified points, which can be computed in one of two ways, the second involving the local Kummer symbols (Proposition~\ref{P:G-primitive elements valuation vector}), and tying everything together.

\section{Preliminary results}
\label{sec:preliminaries}
% (fold)
%%%%%%%%%%%%%%%%%%%%%%%%%%%%%%%%%%%%%%%%%%%%%%%%%%%%%%%%%%%%%%%%%%%%%%%%%%%%%%%%%%%
%%%%%%%%%%%%%%%%%%%%%%%%%%%%%%%%%%%%%%%%%%%%%%%%%%%%%%%%%%%%%%%%%%%%%%%%%%%%%%%%%%%
%%%%%%%%%%%%%%%%%%%%%%%%%%%%%%%%%%%%%%%%%%%%%%%%%%%%%%%%%%%%%%%%%%%%%%%%%%%%%%%%%%%

%%%%%%%%%%%%%%%%%%%%%%%%%%%%%%%%%%%%%%%%%%%%%%%%%%%%%%%%%%%%%%%%%%%%%%%%%%%%%%%%%%%
%%%%%%%%%%%%%%%%%%%%%%%%%%%%%%%%%%%%%%%%%%%%%%%%%%%%%%%%%%%%%%%%%%%%%%%%%%%%%%%%%%%
%%%%%%%%%%%%%%%%%%%%%%%%%%%%%%%%%%%%%%%%%%%%%%%%%%%%%%%%%%%%%%%%%%%%%%%%%%%%%%%%%%%
\subsection{Geometric adeles}
\label{subsec:adeles}
% (fold)
%%%%%%%%%%%%%%%%%%%%%%%%%%%%%%%%%%%%%%%%%%%%%%%%%%%%%%%%%%%%%%%%%%%%%%%%%%%%%%%%%%%
%%%%%%%%%%%%%%%%%%%%%%%%%%%%%%%%%%%%%%%%%%%%%%%%%%%%%%%%%%%%%%%%%%%%%%%%%%%%%%%%%%%
%%%%%%%%%%%%%%%%%%%%%%%%%%%%%%%%%%%%%%%%%%%%%%%%%%%%%%%%%%%%%%%%%%%%%%%%%%%%%%%%%%%

We give a brief summary of the construction of the adeles of a curve. Although this  is analogous to the case of global fields given in~\cite{CaFr}, there are some differences in the geometric case, especially regarding the topology. See~\cite[\S2]{collectanea} for more details.

Let $X$ be a projective, irreducible, non-singular curve over an algebraically closed field $\k$. Let $\Sigma$ be the function field of $X$. We fix the following notation:
\begin{itemize}

% \item $\Sigma$ is the function field of $X$.

\item When we write $x \in X$ it will be implicitly assumed that $x$ is a closed point of $X$, corresponding to the valuation $\upsilon_{x}$ on the function field $\Sigma$. Denote by $\oo_{X,x}$, or simply $\oo_{x}$, the valuation ring at $x$. Since $\k$ is algebraically closed, the closed points are in one-to-one correspondence with the discrete valuations on $\Sigma$.

\item Let $A_{x}$ be the completion of $\oo_{x}$ with respect to $\upsilon_{x}$, which will also denote the extended valuation. Let $\m_{x}$ be the maximal ideal of $A_{x}$, $K_{x}$ its quotient field. Since $\k$ is algebraically closed, the residue field is $A_{x}/\m_{x} = \k$. 

\end{itemize}
A choice of uniformizing parameter $z_{x}$ at $x$ determines the following isomorphisms:% (\cite{Cohen}):
\begin{equation}
\label{E:Ax Kx pws}
	A_{x} \simeq \k\fpws{z_{x}},
	\quad
	K_{x} \simeq \k\laurent{z_{x}},
	\quad
	A_{x}^{*} \simeq \k^{*} \times (1 + \m_{x}).	
\end{equation}

The ring of \emph{adeles} $\A_{X}$ of $\Sigma/\k$ is the subring of $\prod_{x} K_{x}$ given by the restricted direct product with respect to the subrings $A_{x}$, over the closed points of $X$,
\begin{equation*}
\label{E:adele}
	       \A_{X} 
	\eqdef \prodprime_{x \in X} (K_{x},A_{x})
	 = \left\{ (\alpha_{x})_{x\in X} \suchthat \alpha_{x}\in A_{x}  \mbox{ for almost all } x\in X \right\}
\end{equation*}
where ``almost all'' means ``for all but finitely many''. It is equipped with the restricted product topology of the valuation topology on each factor, which we recall is \emph{not} induced by the direct product topology. 
% From now on, we shall always consider $\A_{X}$ as a topological ring and any function field $\Sigma$ will be considered as a topological field with the topology induced by its adele ring; this is in fact the discrete topology, just as for global fields.

% For an adele $\alpha = (\alpha_{x})_{x\in X}$ and a (closed) point $x \in X$ with $\alpha_{x}\in A_{x}$, $\alpha(x)$ will denote the image of its $x$-component $\alpha_{x}$ in the residue field $\k(x)$.

$\A_{X}$ arises as a direct limit as follows: denote by $F$ a finite subset of (closed) points of $X$. Consider
\begin{equation}
\label{E:AXS}
		   \A_{X,F}
	\eqdef \prod_{x \in F} K_{x}
	\times \prod_{x \in X \setminus F} A_{x}.
\end{equation}
% This ring carries a linear topology; namely, the one generated by the neighborhood basis at $0$ consisting of the sets $\prod_{x} \m_{x}^{n_{x}}$, where $(n_{x})_{x\in X}$ runs over collections of non-negative integers $n_{x}$ such that $n_{x}=0$ for almost all $x$.
Given finite subsets $F_1,F_2$ with $F_1\subseteq F_2$, we have an inclusion $\A_{X,F_1}\hookrightarrow \A_{X,F_2}$. Fix a finite subset $F_{0}$ and consider sets containing $F_{0}$. Then the direct limit with respect to these inclusions is isomorphic to $\A_{X}$:
\begin{equation}
\label{E:AX limit}
	\A_{X} \simeq \varinjlim_{F \supseteq F_{0}} \A_{X,F},
\end{equation}
as linearly topological rings. In particular, this does not depend on the choice of $F_{0}$.

% There is an alternative way to generate this topology in the adele ring. Let $\A_{X}^+$ denote the subring of $\A_{X}$ given by the usual direct product
% \[
% 	\A_{X}^+ \eqdef \prod_{x \in X} A_{x},
% \]
% i.e. where the word ``almost'' is dropped in the previous definition. Recalling the general notion of integrality for ring extensions, it is easy to check that $\A_{X}^{+}$ is integrally closed in $\A_{X}$, highlighting the interaction of the algebraic notion of integral closure in commutative rings, with the adelic topology. We will come back to this theme in \S\badref{\ref{subsec:topology}}.
%
%
% For two $\k$-vector spaces  $U,V\subseteq \A_{X}$, we say that they are commensurable and write $U\sim V$ iff
% \[
% 	\dim_{\k} (U+V) / (U\cap V)  <  \infty.
% \]
% Then, the topology of $\A_{X}$ coincides with the topology having  $\{U\subseteq \A_{X}   \suchthat   U \sim \A_{X}^+\}$ as neighborhood basis at $0$. It is worth noticing that $\A_{X}^{+} \subseteq \A_{X,F}$ for all $F$.

Finally, the \emph{idele group} $\I_{X}$ is the group $\A_{X}^{*}$ of invertible elements of $\A_{X}$ endowed with the initial topology of the map
\[
	\I_{X}  \longrightarrow   \A_{X} \times \A_{X}
\]
that sends $\alpha$ to $(\alpha, \alpha^{-1})$. Observe that it is the restricted product of $K_{x}^{*}$ with respect to the unit groups $A_{x}^{*}$.

% and that $\I_X$ can be constructed as the direct limit $\varinjlim \I_{X,F}$ where $\I_{X,F} \eqdef \A^{*}_{X,F}$.

% For a more detailed exposition of the above results, especially regarding the equivalence of the adelic (restricted direct product) and the topology of commensurability, see~\cite[\S2.2]{collectanea}.

% (end)
% end subsection adeles

%%%%%%%%%%%%%%%%%%%%%%%%%%%%%%%%%%%%%%%%%%%%%%%%%%%%%%%%%%%%%%%%%%%%%%%%%%%%%%%%%%%
%%%%%%%%%%%%%%%%%%%%%%%%%%%%%%%%%%%%%%%%%%%%%%%%%%%%%%%%%%%%%%%%%%%%%%%%%%%%%%%%%%%
\subsection{Adelic algebras}
\label{subsec:definition of adelic algebra}
%%%%%%%%%%%%%%%%%%%%%%%%%%%%%%%%%%%%%%%%%%%%%%%%%%%%%%%%%%%%%%%%%%%%%%%%%%%%%%%%%%%
%%%%%%%%%%%%%%%%%%%%%%%%%%%%%%%%%%%%%%%%%%%%%%%%%%%%%%%%%%%%%%%%%%%%%%%%%%%%%%%%%%%

Having briefly reviewed the basic properties of the adele ring $\A_{X}$, we turn to our main objects of study, which we will refer to simply as adelic algebras.  

% The study of the main properties of these algebras will be carried out in \S\ref{sec:preliminaries}. Their structure will allow us to develop an analog of Kummer theory for $\A_{X}$ in \S\ref{sec:Kummer Extensions of AX}. Concrete instances of these type of algebras arise from geometry (see our discussion in \S\ref{sec:future work}) and algebra (Example~\ref{Ex:Heisenberg}).

For the remainder of the paper, $n$ will denote a fixed integer prime to $\chr \k$. When $n$ is assumed prime, it will be denoted by $p$.

\begin{definition}
\label{D:paramvec}
Let $n>1$ be prime to $\chr \k$. A \emph{parameter vector} is an adele $\t = (t_{x})_{x\in X} \in \A_{X} $ such that:
\begin{enumerate}

\item  $t_{x}\neq 0$ for all $x\in X$.

\item The \emph{ramification locus} of $\t$, defined as
\[
    \begin{aligned}
	        \Ram(\t)
	&\eqdef \left\{x\in X \suchthat (n,\upsilon_{x}(t_{x}))\neq n\right\}
	 =      \left\{x\in X \suchthat \upsilon_{x}(t_{x}) \not\equiv 0 \bmod n \right\},
    \end{aligned}
\]
is finite. 

\end{enumerate}
Given a parameter vector $\t$ and $x\in X$, its \emph{ramification index} at $x$ is
\[
	e_{x} \eqdef \frac{n}{(n, \upsilon_{x}(t_{x}))}.
\]
The integer vector $\mathfrak{e}=(e_{x})$ will be called the \emph{ramification profile} of $\t$. Observe that $\Ram(\t) = \{ x \in X \suchthat e_{x} > 1 \}$.
\end{definition}

Note that any \emph{idele} $\t\in\I_X$ serves as a parameter vector. As we shall see below, there are good reasons to restrict ourselves to considering \emph{only} ideles as parameter vectors.

\begin{definition}
\label{D:algebra AXt}
The rank-$n$ adelic algebra over $\A_{X}$ associated to a parameter vector $\t $ is defined as
\begin{equation}\label{E:AXt}
	\AXn{\t} \eqdef \A_{X}[T] / (T^n - \t),
\end{equation}
where $T$ is an indeterminate.
\end{definition}

$\AXn{\t}$ may also be expressed as a direct limit. In order to see this, we need to consider the $K_{x}$-algebra
\begin{equation}
\label{E:Kx}
		\Kxn{t_{x}} \eqdef  K_{x}[T]/(T^n-t_{x}) = \AXn{\t} \otimes_{A_{x}} K_{x}
\end{equation}
for $x \in X$. We also need the analog of the local ring of integers $A_{x}$. For $\upsilon_{x}(t_{x})\geq 0$, which happens at almost all points since any parameter vector $\t$ is in particular an adele by definition, we may consider the subring $\Axn{t_{x}}$ of $\Kxn{t_{x}}$ defined by
\[
	\Axn{t_{x}} \eqdef A_{x}[T]/(T^n-t_{x}).
\]
Choosing $F_0= \{x\in X \suchthat \upsilon_{x}(t_{x})<0\}$ and imitating~\eqref{E:AX limit}, we define, for a finite subset $F \supseteq F_{0}$ of closed points,
\[
		   \AXSn{\t}{F}
	\eqdef \prod_{x \in F} \Kxn{t_{x}}
	\times \prod_{x \in X \setminus F} \Axn{t_{x}}
	\simeq \A_{X,F}[T]/(T^n-\t)
\]
and, taking the direct limit with respect to the inclusion maps, we have shown the following result.

\begin{proposition}
\label{P:AXt restricted product Axtx}
The adelic algebra $\AXn{\t}$ is the restricted direct product of $\Kxn{t_{x}}$ with respect to the subrings $\Axn{t_{x}}$, i.e.
\begin{equation}
\label{E:AXt restricted product Axtx}
\begin{aligned} 
    \AXn{\t} \simeq \varinjlim_{F \supseteq F_{0}} \AXSn{\t}{F}
	         = \prodprime_{x \in X} (\Kxn{t_{x}}, \Axn{t_{x}}).
\end{aligned}
\end{equation}
\end{proposition}

For the time being, we make no mention of topological structure, although this will be relevant in future work.

The reader familiar with the Kummer theory of ring extensions (e.g. as in~\cite{Borevich} or~\cite{Greither}) will recognize that~\eqref{E:AXt} is a candidate for being a Kummer extension of the adele ring $\A_{X}$, although there are several ingredients missing, the most notable of which is perhaps the absence of a group action. Since adelic algebras have infinitely many idempotents and zero divisors, many general references (e.g.~\cite{Ford, Greither, Jan, Nagahara}) which deal with irreducibility, separability and Galois extensions of commutative rings, cannot be directly applied to their study. However, before dealing with this, other basic properties need to be addressed. 

For example, the requirement $t_{x}\neq 0$ for all $x \in X$ in Definition~\ref{D:paramvec} ensures that $\AXn{\t}$ is reduced. We next turn to another important property, namely its separability.

% Regarding separability, one has for example a theorem of Nagahara~\cite[Theorem 2.3]{Nagahara}, relating the separability of a polynomial over $\A_{X}$ with the separability of its projections over the residue fields of $\A_{X}$. However, because of the rather intricate structure of the maximal spectrum of $\A_{X}$ (see~\cite{SerranoHolgado}), this approach unduly complicates matters and will be seen to be unnecessary in the case of the adeles.

% In fact we
% are able to arrive at a characterization only in terms of the closed points of $X$, as part of our ``local-global'' approach to this and other matters.

% , or equivalently, in terms of the projections at each point $x \in X$, which, by definition, have all their coefficients in $A_{x}$ at almost every $x$.

We refer the reader to~\cite{Ford} as a general reference for the theory of separable algebras over commutative rings. Given a commutative ring $R$, a commutative $R$-algebra $A$ is called \emph{separable} over $R$ iff $A$ is projective as an $A\otimes_R A$-module. Furthermore, if $P(T)\in R[T]$ is a monic polynomial, we say that $P(T)$ is separable if $R[T]/(P(T))$ is a separable $R$-algebra. We have the following characterization.

\begin{proposition}\label{P:Axt-separable}
$\AXn{\t}$ is separable as an $\A_{X}$-algebra if and only if $\t\in \I_X$.
\end{proposition}

\begin{proof}
It is straightforward to check that the ideal generated by $P(T)=T^n-\t$ and $P'(T)$ contains $1$ if and only if $\t \in \I_{X}$, from which the conclusion follows by~\cite[Proposition 4.6.1]{Ford}.
% A monic polynomial $P(T) \in R[T]$ over a commutative ring $R$ is separable if and only if the ideal generated by $P$ and its formal derivative $P'$ contains $1$, i.e. $(P(T),P'(T))=R[T]$.		
% It suffices to recall the definition~\eqref{E:AXt} of $\AXn{\t}$ and to apply Proposition~\ref{P:caracterizationseparable} to $P(T)=T^n-\t$: if $\t\in\I_X$, it is clear that $(T^n-\t,nT^{n-1})=1$. Conversely, if the latter holds, there are $A(T),B(T)\in\A_{X}[T]$ with
% \[
% 		A(T)(T^n-\t)+B(T)nT^{n-1}=1.
% \]
% Comparing constant terms one obtains $-a_0\t = 1$, hence  $\t\in\I_X$.
\end{proof}

% For a polynomial $P(T) \in \A_{X}[T]$ and $x \in X$, we will denote by $P(T)_{x} \in K_{x}[T]$ the polynomial whose coefficients are the projections onto $K_{x}$ of the coefficients of $P(T)$.
%
%
% \begin{definition}\label{D:pointwiseseparability}
% Given a monic polynomial $P(T)\in\A_{X}[T]$, $P(T)$ is said to be pointwise separable if  $P(T)_{x}$ is separable over $K_{x}$ for every $x\in X$.
% \end{definition}
%
%
% \begin{proposition}\label{P:sep-point-sep}
% A monic polynomial $P(T)\in\A_{X}[T]$ is separable if and only if is pointwise separable and, for almost every $x\in X$, $P(T)_{x}$ is separable as a polynomial in $A_{x}[T]$.
% \end{proposition}
%
% \begin{proof}
% This is a matter of noting that a Bézout identity $A(T)P(T)+B(T)P'(T)=1$ in $\A_{X}[T]$ implies a corresponding Bézout identity $A(T)_{x} P(T)_{x} + B(T)_{x} P'(T)_{x} = 1$ in $K_{x}[T]$ at each $x \in X$, which is actually an identity in $A_{x}[T]$ at almost every $x$.
% %
% Conversely, given such identities at each $x$, since $K_{x}$ is a field, they can be chosen so that $\deg A(T)_{x}$ and $\deg B(T)_{x}$ are uniformly bounded by $\deg P(T)_{x} = \deg P(T)$, since $P(T)$ is monic. Hence the pointwise Bézout identities can be combined to give an identity in $\A_{X}[T]$.
% \end{proof}

\begin{example}
\label{Ex:Heisenberg}
Let us provide another example of how $\A_{X}$-algebras arise. Let $\t$ be an idele and
$\Phi: \A_{X}^{\oplus n} \to  \A_{X}^{\oplus n}$ the endomorphism of $\A_{X}$-modules given by
\[
\Phi \eqdef 
     \begin{pmatrix}
	 0 & & \ldots & \t \\
	 1 & \ddots & & 0 \\
	 0 & 1 & & \vdots \\
	 0 & \dots & 1 & 0
	 \end{pmatrix}.
\]
Denote by  $\A_{X}[\Phi]$ the $\A_{X}$-algebra of $\End(\A_{X}^{\oplus n})$ generated by $\Phi$. Observe that $\A_{X}[\Phi] \simeq \A_{X}[T]/ (C_{\Phi}(T))$ where $C_{\Phi}(T) =T^n-\t \in \A_{X}[T]$ is the characteristic polynomial of $\Phi$ and hence it is of the type in~\eqref{E:AXt}. This construction is related to the notion of Heisenberg algebra, since $\A_{X}[\Phi]$ is a maximal commutative, separable $\A_{X}$-subalgebra of $\End(\A_{X}^{\oplus n})$, with constant rank $n$ at every $x\in X$ (see for example~\cite[\S4,\S5]{AdamsBergvelt} and~\cite{LiMulase}).
\end{example}

% (end)
% subsection definition of adelic algebra end

%%%%%%%%%%%%%%%%%%%%%%%%%%%%%%%%%%%%%%%%%%%%%%%%%%%%%%%%%%%%%%%%%%%%%%%%%%%%%%%%%%%
%%%%%%%%%%%%%%%%%%%%%%%%%%%%%%%%%%%%%%%%%%%%%%%%%%%%%%%%%%%%%%%%%%%%%%%%%%%%%%%%%%%
\subsection{Isomorphisms}
\label{subsec:isomorphisms of adelic algebras}
% (fold)
%%%%%%%%%%%%%%%%%%%%%%%%%%%%%%%%%%%%%%%%%%%%%%%%%%%%%%%%%%%%%%%%%%%%%%%%%%%%%%%%%%%
%%%%%%%%%%%%%%%%%%%%%%%%%%%%%%%%%%%%%%%%%%%%%%%%%%%%%%%%%%%%%%%%%%%%%%%%%%%%%%%%%%%

In this section, we characterize when two adelic algebras are isomorphic. As one expects, this involves studying isomorphisms between their component algebras at each point $x \in X$. Our main result is Theorem~\ref{T:AXt isomorphic iff ram equal and idele ratio}, which can be thought of as the function field case of~\cite[Lemma 3]{Komatsu01} proved by Komatsu for algebraic number fields. The result is part of the so-called equivalence problem, which studies the degree to which an object attached to a field, in this case, its adele ring, determines it up to isomorphism. We will take up this problem again in a future paper~\cite{adeles03}.

 % We begin with a series of lemmas.

First, we note that for $t \in K_{x}^{*}$, $\Kxn{t} = K_{x}[T]/(T^{n} - t)$ is an étale $K_{x}$-algebra, isomorphic to a product of copies of a cyclic Kummer extension of $K_{x}$. The technical details are given in the following lemma. For a group $G$ a superscript $G^{n}$ will denote the $n$th powers in $G$, while $n$ copies of a set $S$ will be denoted by $\prod^{n} S$. Since we assume that $\k$ is algebraically closed, by~\eqref{E:Ax Kx pws} we have
\begin{equation}
\label{E:n th powers valuation 0 mod n}
	K_{x}^{*n} = \{a \in K_{x}^{*}: \upsilon_{x}(a) \equiv 0 \bmod n\},
\end{equation}
a fact that we will use extensively.

\begin{lemma}
\label{L:structureKxntx}
Let $t \in K_{x}^{*}$, with $m = (n,\upsilon_{x}(t))$ and $e = n/m$. Then $t \in K_{x}^{*m}$, and if we fix an $m$th root $\tau = t^{1/m} \in K_{x}^{*}$ and an $e$th root $\tau^{1/e}$ in some extension field of $K_{x}$, then $K_{x}(\tau^{1/e})$ is a cyclic Kummer extension of degree $e$. Fixing a primitive $m$th root of unity $\xi \in K_{x}$, we obtain an isomorphism
\begin{equation}
\label{E:isomorphism to copies of Kummer}
\begin{split}
	&\psi_{x} : \Kxn{t} = K_{x}[T]/(T^{n} - t) \xrightarrow{\sim} \prod^{m} K_{x}(\tau^{1/e})
	\\
	&\psi_{x}(P(T)) \eqdef (P(\tau^{1/e}),P(\xi\tau^{1/e}),\dots,P(\xi^{m-1}\tau^{1/e}))
\end{split}
\end{equation}
mapping $K_{x}$ into the diagonal. In particular, we note the following two special cases (the only ones if $n$ is prime):
\begin{enumerate}

\item (Unramified case) If $\upsilon_{x}(t) \equiv 0 \bmod n$, i.e. $e = 1$, then 
% $t \in K_{x}^{*n}$, and choosing $\tau = t^{1/n} \in K_{x}^{*}$ and a primitive $n$th root of unity $\xi \in K_{x}$, we obtain an isomorphism
%%
\begin{equation}
\label{E:isomorphism to copies of Kx}
\begin{split}
	&\psi_{x} : \Kxn{t} = K_{x}[T]/(T^{n} - t) \simeq \prod^{n} K_{x}
	\\
	&\psi_{x}(P(T)) = (P(\tau),P(\xi\tau),\dots,P(\xi^{n-1}\tau)).
\end{split}
\end{equation}
%%
% mapping $K_{x}$ to the diagonal.

\item (Totally ramified case) If $(n,\upsilon_{x}(t)) = 1$, i.e. $e = n$, then $\Kxn{t}$ is a field, isomorphic to the cyclic Kummer extension $K_{x}(t^{1/n})$ of degree $n$.

\end{enumerate}
\end{lemma}

\begin{proof}
This follows from Capelli's Theorem, which states that in a field $K$ containing the $n$th roots of unity $\mu_{n}$, given $a \in K^{*}$, the polynomial $T^{n} - a \in K[T]$ is irreducible if and only if $a \notin K^{*p}$ for any prime $p \mid n$ and $a \notin -4 K^{*4}$ if $4 \mid n$.

In our case, since $m \mid \upsilon_{x}(t)$, we have $t \in K_{x}^{*m}$ by~\eqref{E:n th powers valuation 0 mod n}. Choosing $\tau = t^{1/m} \in K^{*}_{x}$ and recalling that $n$ is assumed prime to the characteristic of $\k$, Capelli's Theorem shows that the factorization of $T^{n} - t$ into irreducibles in $K_{x}[T]$ is
\[
	  T^{n} - t
	= \prod_{\xi \in \mu_{m}} (T^{e} - \xi \tau).
\]
Fix a primitive $m$th root of unity $\xi \in K_{x}$. Then $\mu_{m} = \langle \xi^{e} \rangle$ and hence
\[
	K_{x}[T]/(T^{n} - t) \simeq \prod_{i=1}^{m} K_{x}[T]/(T^{e} - \xi^{ei} \tau).
\]
Now let $\tau^{1/e}$ denote any $e$th root in an extension field of $K_{x}$. Then $L = K_{x}(\tau^{1/e})$ is a cyclic Kummer extension of $K_{x}$ of degree $e$ and $L = K_{x}(\xi^{i}\tau^{1/e})$ for any $i$. Clearly $T^{e} - \xi^{ei}\tau$ is the minimal polynomial of $\xi^{i} \tau^{1/e}$ over $K_{x}$, hence the map $P(T) \mapsto P(\xi^{i}\tau^{1/e})$ induces an isomorphism $K_{x}[T]/(T^{e} - \xi^{ei} \tau) \simeq L$, and therefore~\eqref{E:isomorphism to copies of Kummer}
induces an isomorphism $\Kxn{t} = K_{x}[T]/(T^{n} - t) \simeq \prod^{m} L$.
\end{proof}

\begin{corollary}
\label{C:KxtxnIsomex}
For $t_{i,x} \in K_{x}^{*}$, $i=1,2$, the $K_{x}$-algebras $\Kxn{t_{1,x}}$ and $\Kxn{t_{2,x}}$ are isomorphic if and only if $e_{1,x} = e_{2,x}$.
\end{corollary}

% We highlight the following two special cases, which, when $n$ is prime, are in fact the only ones.

% \begin{corollary}
% \label{C:structure Kxnt e=1}
% Let $t \in K_{x}^{*}$. Then
% %%
% \begin{enumerate}
%
% \item (Unramified case) If $\upsilon_{x}(t) \equiv 0 \bmod n$, i.e. $e = 1$, then $t \in K_{x}^{*n}$, and choosing $\tau = t^{1/n} \in K_{x}^{*}$ and a primitive $n$th root of unity $\xi \in K_{x}$, we obtain an isomorphism
% %%
% \begin{equation}
% \label{E:isomorphism to copies of Kx}
% \begin{split}
% 	&\psi_{x} : \Kxn{t} = K_{x}[T]/(T^{n} - t) \simeq \prod^{n} K_{x}
% 	\\
% 	&\psi_{x}(P(T)) = (P(\tau),P(\xi\tau),\dots,P(\xi^{n-1}\tau))
% \end{split}
% \end{equation}
% %%
% mapping $K_{x}$ to the diagonal.
%
% \item (Totally ramified case) If $(n,\upsilon_{x}(t)) = 1$, i.e. $e = n$, then $\Kxn{t}$ is a field, isomorphic to the cyclic Kummer extension $K_{x}(t^{1/n})$ of degree $n$.
%
% \end{enumerate}
% %%
% \end{corollary}

% \begin{remark}
% \label{R:parameter vectors}
% In particular, Lemma~\ref{L:structureKxntx} shows that for $n > 1$ and $\t = (t_{x})_{x\in X}\in \A_{X}$, the $\A_{X}$-algebra $\A_{X}[T]/(T^n-\t)$ is reduced if and only if $t_{x}\neq 0$ for all $x \in X$, which is why we have included this condition in the definition of a parameter vector (Definition~\ref{D:paramvec}).
% \end{remark}

The following technical result is probably known, but we have not been able to find a reference in the literature, hence we provide a self-contained proof.

\begin{proposition}
\label{P:valuationsAndIsomorphisms}
Given $t_{i,x} \in K_{x}^{*}$ for $i=1,2$ with $\upsilon_{x}(t_{1,x}),\upsilon_{x}(t_{2,x})\geq 0$ and $\upsilon_{x}(t_{1,x}) \equiv \upsilon_{x}(t_{2,x}) \equiv 0\bmod n$. Then there exists a $K_{x}$-algebra isomorphism $\phi_{x}:\Kxn{t_{1,x}}\rightarrow \Kxn{t_{2,x}}$ sending $\Axn{t_{1,x}}$ to $\Axn{t_{2,x}}$, if and only if $\upsilon_{x}(t_{1,x})= \upsilon_{x}(t_{2,x})$.
\end{proposition}

\begin{proof}
If $\upsilon_{x}(t_{1,x}) = \upsilon_{x}(t_{2,x})$, then $\phi_{x}(T)= \tau T$, where $\tau\in A_{x}^*$ satisfies $\tau^n = t_{1,x}/t_{2,x}$, is a $K_{x}$-algebra isomorphism from $\Kxn{t_{1,x}}$ to $\Kxn{t_{2,x}}$ that sends $\Axn{t_{1,x}}$ to $\Axn{t_{2,x}}$.

For the reverse implication, let  $\phi_{x}:\Kxn{t_{1,x}}\rightarrow \Kxn{t_{2,x}}$ be an isomorphism sending $\Axn{t_{1,x}}$ to $\Axn{t_{2,x}}$. Choose, as in Lemma~\ref{L:structureKxntx}, a primitive $n$th root of unity $\xi$ and for $i=1,2$, an $n$th root $\tau_{i}$ of $t_{i,x}$, and let $\psi_{i,x}$ be the corresponding isomorphisms as in~\eqref{E:isomorphism to copies of Kx}.

Since $\phi_{x}$ maps $\Axn{t_{1,x}}$ to $\Axn{t_{2,x}}$, there is a unique polynomial $P(T) = \sum a_{i} T^{i} \in A_{x}[T]$ of degree less than $n$ such that $\phi_{x}(T) = P(T)$ in $\Axn{t_{2,x}}$. Now, the automorphism $\sigma \eqdef \psi_{2,x} \circ \phi_{x} \circ (\psi_{1,x})^{-1}$ of $\prod^{n} K_{x}$ must be a permutation of the factors, and the relation
\[
		\sigma(\psi_{1,x}(T))= \psi_{2,x}(P(T))
\]
translates, in light of~\eqref{E:isomorphism to copies of Kx}, to
\[
    \begin{aligned}
       \sigma(\tau_{1}, \xi \tau_{1}, \dots, \xi^{n-1} \tau_{1})
	&= (\xi^{s_0}\tau_1, \dots ,\xi^{s_{n-1}}\tau_1)
	\\
	&= (P(\tau_{2}), P(\xi \tau_{2}), \dots, P(\xi^{n-1} \tau_{2})),
    \end{aligned}
\]
for some permutation $\{s_0, \dots ,s_{n-1}\}$ of $\{0, \dots ,n-1\}$. This leads to the system of equations
\[
	a_{0} + a_{1} (\xi^{k} \tau_{2}) + \dots + a_{n-1} (\xi^{k} \tau_{2})^{n-1}
	= \xi^{s_{k}} \tau_{1}
	\quad
	(k = 0,1,\dots,n-1)
\]
and by Cramer's rule,
\[
		a_{k}= \frac{\tau_1}{\tau_2^{k}} c_{k},
\]
where $c_{k}$ is a rational function of $\xi$ and hence lies in $\k$. Since $P(T)$ is nonconstant, $a_{k} \neq 0$ for some $k \neq 0$, therefore
\[
		  \upsilon_{x}(a_{k}^{n}) 
		= \upsilon_{x}(\tau_{1}^{n}) - k \upsilon_{x} (\tau_{2}^{n}) 
		\geq 0
		\implies \upsilon_{x} (t_{x,1}) \geq \upsilon_{x}(t_{x,2})
\]
and the reverse inequality follows by considering $\phi_{x}^{-1}$.
\end{proof}

\begin{lemma}
\label{L:familiesOfIsomorphisms}
The $\A_{X}$-algebra isomorphisms $\phi:\AXn{\t_{1}}\rightarrow \AXn{\t_{2}}$ correspond bijectively to families $\{\phi_{x}:\Kxn{t_{1,x}}\rightarrow \Kxn{t_{2,x}}\}_{x\in X}$ of $K_{x}$-algebra isomorphisms such that $\phi_{x}(\Axn{t_{1,x}})= \Axn{t_{2,x}}$ for almost every $x\in X$.
\end{lemma}

\begin{proof}
Given an $\A_{X}$-algebra isomorphism $\phi$, $x \in X$, and $\lambda \in \Kxn{t_{1,x}}$, define $\phi_{x}(\lambda) = (\phi(\alpha))_{x} \in  \Kxn{t_{2,x}}$, where $\alpha \in \AXn{\t_{1}}$ is any element with $\alpha_{x}= \lambda$. Note that, given $x\in X$, the element of $\A_{X}$ defined by $1\in \Kxn{t_{1,x}}$ and $0 \in \Kxn{t_{1,y}}$ for $y\neq x$ is an idempotent and, thus, the $x$-component of $\phi(\alpha)$ only depends on the $x$-component of $\alpha$.  Hence, $\phi_{x}$ is well-defined.

Note that $\phi_{x}(\Axn{t_{1,x}}) \subseteq \Axn{t_{2,x}}$ for almost every $x\in X$; otherwise there would be infinitely many $x \in X$ and corresponding $\lambda_{x} \in \Axn{t_{1,x}}$ such that $\phi_{x}(\lambda_{x}) \notin \Axn{t_{2,x}}$, allowing us to choose $\alpha \in \AXn{\t_{1}}$ with $\phi(\alpha) \notin \AXn{\t_{2}}$. Since the same reasoning can be applied to $\phi^{-1}$, the desired equality holds at almost every $x$.

Conversely, if $\{\phi_{x}:\Kxn{t_{1,x}}\rightarrow \Kxn{t_{2,x}}\}_{x\in X}$ is now a family as described, the map $\phi:\AXn{\t_{1}}\rightarrow\AXn{\t_{2}}$ defined by $(\phi(\alpha))_{x}= \phi_{x}(\alpha_{x})$ is a well-defined $\A_{X}$-isomorphism.
\end{proof}

\begin{theorem}
\label{T:AXt isomorphic iff ram equal and idele ratio}
Given two parameter vectors $\t_{1}$ and $\t_{2}$, the corresponding $\A_{X}$-algebras $\AXn{\t_{1}}$ and $\AXn{\t_{2}}$ are isomorphic if and only if $\mathfrak{e}_{1}= \mathfrak{e}_{2}$ and $\t_{1}  = u \cdot  \t_{2}$ for some $u\in \I_X$. In particular, if $\t_{1},\t_{2}\in\I_X$, then $\AXn{\t_{1}}\simeq\AXn{\t_{2}}$ if and only if $\mathfrak{e}_{1}= \mathfrak{e}_{2}$.
\end{theorem}

\begin{proof}
If $\AXn{\t_{1}}$ and $\AXn{\t_{2}}$ are isomorphic, by Corollary~\ref{C:KxtxnIsomex} and Lemma~\ref{L:familiesOfIsomorphisms}, we obtain $\mathfrak{e}_{1}= \mathfrak{e}_{2}$. The idele $u\in\I_X$ appearing in the statement can be defined by $u_{x} = t_{1,x} / t_{2,x}$; this has valuation $0$ at almost every point by Lemma~\ref{L:familiesOfIsomorphisms} and Proposition~\ref{P:valuationsAndIsomorphisms}. 

Conversely, at almost every $x$, by definition, we have $\upsilon_{x}(t_{1,x}),  \upsilon_{x}(t_{2,x}) \geq 0$, $\upsilon_{x}(t_{1,x}) \equiv \upsilon_{x}(t_{2,x}) \equiv 0 \bmod n$, and $\upsilon_{x}(u_{x}) = 0$ since $u \in \I_{X}$. At such points, $\upsilon_{x}(t_{1,x}) = \upsilon_{x}(t_{2,x})$, hence by Proposition~\ref{P:valuationsAndIsomorphisms} we have a $K_{x}$-algebra isomorphism $\phi_{x}:\Kxn{t_{1,x}}\rightarrow \Kxn{t_{2,x}}$ such that $\phi_{x}(\Axn{t_{1,x}})= \Axn{t_{2,x}}$. At any other points $x$, of which there are only a finite number, since $e_{1,x} = e_{2,x}$, by Corollary~\ref{C:KxtxnIsomex}, there is a $K_{x}$-algebra isomorphism $\phi_{x}:\Kxn{t_{1,x}}\rightarrow \Kxn{t_{2,x}}$. By Lemma~\ref{L:familiesOfIsomorphisms}, the family $\{\phi_{x}\}_{x}$ defines an $\A_{X}$-algebra isomorphism between $\AXn{\t_{1}}$ and $\AXn{\t_{2}}$.
\end{proof}

\section{\texorpdfstring{Kummer Extensions of $\A_{X}$}{Kummer Extensions of the adele ring}}
\label{sec:Kummer Extensions of AX}
% (fold)
%%%%%%%%%%%%%%%%%%%%%%%%%%%%%%%%%%%%%%%%%%%%%%%%%%%%%%%%%%%%%%%%%%%%%%%%%%%%%%%%%%%
%%%%%%%%%%%%%%%%%%%%%%%%%%%%%%%%%%%%%%%%%%%%%%%%%%%%%%%%%%%%%%%%%%%%%%%%%%%%%%%%%%%
%%%%%%%%%%%%%%%%%%%%%%%%%%%%%%%%%%%%%%%%%%%%%%%%%%%%%%%%%%%%%%%%%%%%%%%%%%%%%%%%%%%

A general Galois theory for arbitrary commutative rings was developed in the 1960s, starting with the seminal work of Chase, Harrison and Rosenberg~\cite{ChaseHarrisonRosenberg}. It can be applied to the case we are interested in here, namely, an $\A_{X}$-algebra $S$ with a group action via $\A_{X}$-automorphisms. A distinctive feature of the Galois theory of rings is that in general there is no unique ``Galois group'', but rather, a fixed extension may be Galois under the action of many different groups.

In this section we will develop the analog of Kummer theory for the adele ring of the function field of an algebraic curve as the first step towards the general abelian Galois theory. Specifically, we will be considering ring homomorphisms $\A_{X} \to S$ which are Galois under a fixed action of a cyclic group of prime order. In addition to the reference mentioned above, we also refer to~\cite{Borevich,Ford,Greither}.

Based on the existence of primitive elements for this type of ring extensions, we define invariants which explicitly classify the $p$-cyclic Galois extensions of $\A_{X}$.  We also give a complete determination of the admissible Galois structures for a given fixed extension of $\A_{X}$.

As can reasonably be expected for the ring $\A_{X}$, there are relations between the global Galois property and the corresponding ones at each point $x \in X$.

%%%%%%%%%%%%%%%%%%%%%%%%%%%%%%%%%%%%%%%%%%%%%%%%%%%%%%%%%%%%%%%%%%%%%%%%%%%%%%%%%%%
%%%%%%%%%%%%%%%%%%%%%%%%%%%%%%%%%%%%%%%%%%%%%%%%%%%%%%%%%%%%%%%%%%%%%%%%%%%%%%%%%%%
%%%%%%%%%%%%%%%%%%%%%%%%%%%%%%%%%%%%%%%%%%%%%%%%%%%%%%%%%%%%%%%%%%%%%%%%%%%%%%%%%%%
\subsection{\texorpdfstring{Structure of $p$-cyclic extensions}{Structure of p-cyclic extensions}}
\label{subsec:structure of p cyclic extensions}
% (fold)
%%%%%%%%%%%%%%%%%%%%%%%%%%%%%%%%%%%%%%%%%%%%%%%%%%%%%%%%%%%%%%%%%%%%%%%%%%%%%%%%%%%
%%%%%%%%%%%%%%%%%%%%%%%%%%%%%%%%%%%%%%%%%%%%%%%%%%%%%%%%%%%%%%%%%%%%%%%%%%%%%%%%%%%
%%%%%%%%%%%%%%%%%%%%%%%%%%%%%%%%%%%%%%%%%%%%%%%%%%%%%%%%%%%%%%%%%%%%%%%%%%%%%%%%%%%

In this section we give a brief overview of the Galois theory of commutative rings, referring the reader to the above mentioned sources. The main result is Theorem~\ref{T:galois local global}, which establishes a local-global criterion for when a ring extension of $\A_{X}$ is Galois for a $p$-cyclic group $G$, in terms of the closed points on the curve $X$. In addition it also yields the $G$-module structure of such an extension. The proof relies on the triviality of the Picard group of the ring of adeles, which is shown in Theorem~\ref{T:Pic adeles}.

There are many ways to define a Galois extension of rings. A list of equivalent properties is given in~\cite[Theorem 12.2.9]{Ford}. We choose the following one, which is well suited to our point of view.
\begin{definition}
\label{D:strongly distinct}
Two morphisms of commutative rings $f,g: S \to T$ are called \emph{strongly distinct} if for every nonzero idempotent $e\in T$  there exist $s\in S$ such that $f(s)e \neq g(s)e$.
\end{definition}

\begin{definition}[Galois extension of rings]
\label{D:Galois extension of rings}
A Galois extension of a commutative ring $R$ consists of a pair $(S,G)$, where $S$ is a commutative ring extension of $R$ and $G$ is a finite group such that:
\begin{itemize}

\item $G$ acts faithfully on $S$ by $R$-algebra automorphisms.

\item $S^{G} = R$.

\item $S$ is a separable $R$-algebra.

\item The elements of $G$ are pairwise strongly distinct.

\end{itemize}
In this case we say that $S$ is a Galois ring extension of $R$ with Galois group $G$, or simply a $G$-Galois extension of $R$.
\end{definition}

It is easily checked that a Galois extension of fields in the classical sense is also a Galois extension of rings under this general definition.

An important special case, given in the following example from~\cite{Greither}, generalizes the concept of a radical Kummer extension of fields. We have modified it slightly so that the group action is via characters. In addition, we need to restrict to base rings $R$ of the following type.

\begin{definition}[Kummerian ring]
\label{D:Kummerian ring}
Let $n$ be a natural number. A commutative ring $R$ is $n$-\emph{Kummerian} if $n$ is prime to $\chr(R)$ and its unit group $R^{*}$ contains a distinguished $n$-cyclic subgroup $\mu_{n}$.
\end{definition}

Note that if $K$ is a field with $\chr(K)$ prime to $n$ and containing the $n$-th roots of unity $\mu_{n} = \mu_{n}(K^{*})$, then $K$ is $n$-Kummerian and this subgroup is the \emph{only} possible choice.

For our purposes, since we have fixed the algebraically closed field $\k$, the choice of $\mu_{n} \eqdef \mu_{n}(\k^{*})$, the group of $n$th roots of unity in $\k$, induces on any $\k$-algebra $R$ of characteristic prime to $n$ the structure of an $n$-Kummerian ring, which is compatible with $\k$-algebra morphisms.

In particular this will be the case for $R = \A_{X}$, where we have $\k \subseteq \Sigma$ embedded diagonally in $\A_{X}$ and copies of $\k$ in each completion $K_{x}$. This example highlights how there may be infinitely many choices of subgroups $\mu_{n} \subseteq R^{*}$, hence the necessity of specifying one beforehand.

\begin{example}
\label{EX:G-chi Kummer}
Let $R$ be an $n$-Kummerian ring, $u \in R^{*}$ a unit, $G$ an $n$-cyclic group, and $\chi : G \to \mu_{n} \subseteq R^{*}$ a character of order $n$. The extension
\[
	\gkummer{R}{u} \eqdef R[T]/(T^{n} - u)
\]
may be endowed with an action of $G$ via $\chi$ by defining
\begin{equation}
\label{E:chi-action on T}
		g(T) \eqdef \chi(g) T.
\end{equation}
As shown in~\cite[p.20]{Greither}, $\gkummer{R}{u}/R$ is a $G$-Galois extension of rings.
\end{example}

In particular, this applies to an adelic algebra $\AXn[n]{\t}$, where $\t \in \I_{X}$ is an idele, giving it the structure of an $n$-cyclic ring extension of $\A_{X}$.

Since the Galois extensions considered in Example~\ref{EX:G-chi Kummer} play a central role in the Galois theory of $\A_{X}$, it is useful to give them a name.

\begin{definition}[$(G, \chi)$-Kummer extensions]
\label{D:G-chi-Kummer extension}
For a fixed $n$-Kummerian base ring $R$, a \emph{$(G,\chi)$-Kummer extension of $R$} is a triple $(\gkummer{R}{u},G,\chi)$, where $\gkummer{R}{u} \eqdef R[T]/(T^{n} - u)$ with $u \in R^{*}$, and $G$ is an $n$-cyclic group which acts on $\gkummer{R}{u}$ via the character $\chi : G \to \mu_{n} \subseteq R^{*}$ of order $n$ as in~\eqref{E:chi-action on T}.
\end{definition}

We review some of the basic facts regarding group actions on modules over a ring $R$. Let $R$ be an $n$-Kummerian ring with distinguished subgroup $\mu_{n}$, and $S$ an $R$-module. Suppose $G$ is a finite abelian group of order $n$ acting on $S$ via $R$-module automorphisms. Its dual group $\widehat{G}$ will be identified with $\Hom(G,\mu_{n})$ and its elements referred to simply as characters of $G$. We may consider the decomposition of $S$ with respect to the action of $\widehat{G}$, namely, for $\chi \in \widehat{G}$, we define the $\chi$-eigenspace (or isotypical component)
\begin{equation*}
\label{E:chi eigenspace definition}
	S^{\chi} \eqdef \{ \alpha \in S \suchthat g(\alpha) = \chi(g) \alpha \ \forall g \in G \}.
\end{equation*}
Projection onto the $\chi$-eigenspace is given by $\alpha_{\chi} = e_{\chi} \alpha$, where $e_{\chi}$ is the corresponding idempotent in the group algebra,
\begin{equation*}
\label{E:e-chi idempotent}
	e_{\chi} \eqdef \frac{1}{|G|} \sum_{g \in G} \chi(g^{-1}) g \in R[G].
\end{equation*}
We then have the decomposition
\begin{equation}
\label{E:eigenspace decomposition}
	S = \bigoplus_{\chi\in\widehat{G}} S^{\chi}.
\end{equation}

For simplicity, we will consider from now on that the rank $n$ is a prime number $p$, different from the characteristic of $\k$.

If $G$ is cyclic of prime order and the action is nontrivial, then each $R$-module $S^{\chi}$ is nontrivial.

Going back to the case of base ring $R = \A_{X}$, we explore some of the relations between the properties of an extension $S$ of $\A_{X}$ and the corresponding properties of its localizations at the closed points of $X$.

\begin{proposition}
\label{P:eigenspace commutes with tensor product}	
Let $S$ be an $\A_{X}$-module. Suppose that $G$ is a $p$-cyclic group acting on $S$ via $\A_{X}$-module automorphisms, with $p$ prime and different from $\chr{\k}$. 
For any point $x \in X$ and any character $\chi : G \to \mu_{p} \subseteq \k^{*}$, we have
\begin{equation}
\label{E:eigenspace commutes with tensor product}
	(S \otimes_{\A_{X}} K_{x})^{\chi} = S^{\chi} \otimes_{\A_{X}} K_{x}.
\end{equation}
\end{proposition}

\begin{proof}
Keeping in mind that $G$ acts on $S$ via $\A_{X}$-automorphisms and trivially on $K_{x}$, and that the maps $S \to S \otimes_{\A_{X}} K_{x}$ are surjective, one has the following commutative diagram:
\begin{equation*}
\begin{tikzcd}
 	S \arrow[d, two heads,  "e_{\chi}"']  \arrow[r, two heads] & S \otimes_{\A_{X}} K_{x} \arrow[d, two heads, "e_{\chi}"] 
 	\\
    S^{\chi} \arrow[r]  & (S \otimes_{\A_{X}} K_{x})^{\chi}
\end{tikzcd}
\end{equation*} 
hence the bottom arrow is also a surjection. In addition, this map factors via the inclusion
\[
	S^{\chi} \xrightarrow{} S^{\chi} \otimes_{\A_{X}} K_{x} \hookrightarrow (S \otimes_{\A_{X}} K_{x})^{\chi}
\]
from which~\eqref{E:eigenspace commutes with tensor product} follows.
\end{proof}

Henceforward we will use the notation $S_{x} \eqdef S \otimes_{\A_{X}} K_{x}$.

\begin{lemma}
\label{L:S injects prod Sx}
Let $P$ be a finitely generated projective $\A_{X}$-module. Then $P$ injects into the direct product $\prod_{x} P_{x}$.
\end{lemma}

\begin{proof}
The result is clear when $P = \A_{X}$, and hence for free $\A_{X}$-modules of finite rank. The general result follows since a finitely generated projective module is a direct summand of a finitely generated free module.
% For each finite subset $F \subseteq X$, denote $I_{F} = \ker(\A_{X} \to \prod_{x \in F} K_{x})$. We have the exact sequence
% \[
% 	0 \xrightarrow{} I_{F} \xrightarrow{} \A_{X} \xrightarrow{} \A_{X}/I_{F} \xrightarrow{} 0.
% \]
% Tensoring by $S$, which is a finite flat $\A_{X}$-module, we have
% \[
% 	0 \xrightarrow{} I_{F} S \xrightarrow{} S \xrightarrow{} S /I_{F} S = \prod_{x \in F} S_{x}
% 	\xrightarrow{} 0,
% \]
% which forms an inverse system of exact sequences. Noting that $I_{F}$ is principal, generated by the idempotent which has a $0$ at the places in $F$ and $1$ elsewhere, we see that
% \[
% 	I_{F} S \cap I_{F'} S = I_{F} I_{F'} S = I_{F \cup F'} S
% \]
% and in particular $\bigcap_{F} I_{F} S = 0$, hence taking the limit, we obtain
% \[
% 	0 \xrightarrow{} \varprojlim_{F} I_{F}S = 0
% 	  \xrightarrow{} S
% 	  \xrightarrow{} \varprojlim_{F} \prod_{x \in F} S_{x} = \prod_{x} S_{x}
% \]
% and we conclude.
\end{proof}

In the case of the ring of adeles, we are able to determine the $p$-cyclic Galois extensions of $\A_{X}$ from the triviality of its Picard group. The latter can  be deduced from the very general development in~\cite[\S2]{Cesnavicius} or~\cite{Bhatt}, but we prefer to include a proof that only requires standard commutative algebra, in order to make the exposition self-contained.

\begin{theorem}
\label{T:Pic adeles}	
$\Pic(\A_{X})=0$.
\end{theorem}

\begin{proof}
We begin by showing that $\Pic(\A_{X,F}) = 0$, where the ring $\A_{X,F}$ is defined in~\eqref{E:AXS} for a finite subset $F \subseteq X$. One has
\[
	\Pic( \A_{X,F})=
	\prod_{x\in F} \Pic( K_{x} ) \times \Pic(\prod_{x\in X \setminus F} A_{x}).
\]
Since $K_{x}$ is a field, $\Pic( K_{x} )=0$, and thus it suffices to show that $\Pic(R)=0$ where $R=\prod_{i\in I} R_{i}$ for any set $I$ of DVRs $R_{i}$. Let $L\in \Pic(R)$.  Since $L$ is locally free and $\Spec(R)$ is compact, there is a finite set of elements $f_1,\ldots, f_n\in R$ such that $U_j\eqdef\Spec(R_{f_{j}})= \Spec(R)\setminus(f_{j})_{0}$ is a covering of $\Spec(R)$ and the localization $L_{f_{j}}$ is a rank $1$ free $R_{f_{j}}$-module. Let $L_i\eqdef L\otimes_{R} R_{i} $ which lies in $\Pic(R_{i})$. Since $L$ is finitely generated, then there is a canonical map of $R$-modules
	\begin{equation}
	\label{E:L to L_i}
	L \isomto L\otimes_{R} (\prod R_{i}) \to \prod (L\otimes_{R}  R_{i}) \isomto \prod L_{i}.
	\end{equation}
Noting that $R_i$ is a DVR and that $L_i$ is locally free of rank $1$, it follows that $L_i$ is free of rank $1$ and, thus, $\prod L_i$ is a rank $1$ free $R$-module. Furthermore, when restricted to $U_j$, both $L$ and $\prod L_i$ are free modules, and the map \eqref{E:L to L_i} is an isomorphism. Hence, $L$ is free. 

Now, recall from~\eqref{E:AX limit} that $\A_{X}=\varinjlim\A_{X,F}$. We claim that
\[
	\Pic(\A_{X}) = \varinjlim \Pic( \A_{X,F}).
\]
Let $L\in \Pic(\A_{X}) $. By \cite[Theorem 8.5.2.(ii)]{EGA IV.3}, there exists $F$ and a quasi-coherent $\A_{X,F}$-module of finite presentation $L_F$ such that $L=L_F\otimes_{\A_{X,F}}\A_X$. For $F\subset F'$, set $L_{F'}\eqdef L_F\otimes_{\A_{X,F}}\A_{X,F'}$. Applying \cite[Proposition 8.5.5]{EGA IV.3}, it follows that there exists $F'$ such that $L_{F'}\in \Pic( \A_{X,F'})$ and, thus $L_{F''} \in  \Pic( \A_{X,F''})$  holds also for any $F''$ containing $F'$.  Conversely, given $\{L_F\}\in \varinjlim \Pic( \A_{X,F})$. Then, $L \eqdef L_F\otimes_{\A_{X,F}}\A_X $, which does not depend on $F$, belongs to $\Pic(\A_X)$. 
\end{proof}

\begin{lemma}
\label{L:density points Spec adeles}
The set of closed points $x \in X$ is dense in $\Spec(\A_{X})$ with the Zariski topology.
\end{lemma}

\begin{proof}
We show that for each open subset $U\subseteq \Spec \A_{X}$ there is a closed point $x\in X$ such that the maximal ideal of $\A_{X}$ defined by $x$, $I_{x} \eqdef \ker(\A_{X} \to K_{x})$, lies in $U$. Let $U$ be an open subset defined by the zeroes of an ideal of $\A_{X}$. If $U= \Spec \A_{X}$, it is clear that the set of closed points of $X$ belong to $U$. Assume that $U$ is neither empty nor the whole space. Let $Z$ be the complement of $U$ which, being closed, is the set of zeroes $Z=(I)_{0}$ of an ideal $I$ with $(0)\subsetneq I \subsetneq \A_{X}$.  Thus, there exists $\alpha\in I \setminus\{0 \}$ which is not invertible, and hence there is an $x\in X$ such that $\alpha_x\neq 0$ and $I \not\subset I_{x}$. Therefore $I_{x} \notin Z$, i.e. $ I_{x} \in U$.
\end{proof}

% For a finitely generated algebra $S$ over any commutative ring $R$, it is known that $S$ is separable if and only if its localizations at every maximal ideal are separable (\cite[Theorem 7.1]{Meyer}). The corresponding local-global criterion for the separability of a monic polynomial is given in~\cite[Theorem 2.3]{Nagahara}. In the case of base ring $R = \A_{X}$, Proposition~\ref{P:sep-point-sep} gives a criterion for separability which involves only the closed points $x \in X$, which are a dense part but not all of $\Spec(\A_{X})$ (see~\cite{SerranoHolgado}).

The next theorem establishes an analogous local-global relation which characterizes when an $\A_{X}$-algebra extension of $\A_{X}$ is Galois for a given $p$-cyclic group $G$. In particular we see that the conditions of Definition~\ref{D:Galois extension of rings} may be simplified in the case of $R = \A_{X}$.

\begin{theorem}
\label{T:galois local global}
Let $S$ be a separable $\A_{X}$-algebra which is finitely generated and projective as an $\A_{X}$-module, and $G$ a $p$-cyclic group with $p$ a prime different from $\chr(\k)$, acting on $S$ by $\A_{X}$-automorphisms. Then the following are equivalent:
\begin{enumerate}

\item $S$ is a $G$-Galois ring extension of $\A_{X}$.

\item $S_{x}$ is a $G$-Galois ring extension of $K_{x}$ for all closed points $x \in X$.

\end{enumerate}
If this is the case, then $S^{\chi}$ is a free $\A_{X}$-module of rank $1$ for all characters $\chi : G \to \mu_{p} \subseteq \k^{*}$, and thus $S$ is a free $\A_{X}$-module of rank $p$.
\end{theorem}

\begin{proof}
If $S$ is $G$-Galois over $\A_{X}$, then by the base change property for Galois extensions of commutative rings (Lemma 12.2.7(1) of~\cite{Ford}) we immediately conclude that each $S_{x}$ is $G$-Galois over $K_{x}$.

For the converse, we need to verify the four conditions in Definition~\ref{D:Galois extension of rings}. Since separability is assumed, that leaves the other three.

\begin{itemize}

\item To see that $S^{G} = S^{\chi_{0}} = \A_{X}$ we will in fact show that (2) already implies that each eigenspace $S^{\chi}$ is a free $\A_{X}$-module of rank $1$, and hence by~\eqref{E:eigenspace decomposition} $S$ is free of rank $p$.

Since $S$ is finitely generated and projective and $S^{\chi}$ is a direct summand of $S$ (see~\eqref{E:eigenspace decomposition}), it follows that $S^{\chi}$  is a finite and locally free $\A_{X}$-module (\cite[\href{https://stacks.math.columbia.edu/tag/00NX}{Lemma 00NX(3)}]{stacks-project}). By Theorem~\ref{T:Pic adeles} it will suffice to check that it has rank $1$. 

From~\cite[\href{https://stacks.math.columbia.edu/tag/00NX}{Lemma 00NX(8)}]{stacks-project}, we know that the function mapping $\mathfrak{p}\in \Spec \A_{X}$ to $\dim_{\k(\mathfrak{p})} S^{\chi}_{\k(\mathfrak{p})} $ is locally constant in the Zariski topology. By Lemma~\ref{L:density points Spec adeles}, it is enough to show that $\dim_{K_{x}} S^{\chi}_{x}=1$ for any (closed) point $x \in X$.

Observe first that by base change for separable ring extensions (see~\cite[\S4.3]{Ford}),
$S_x$ is a separable $K_{x}$-algebra since $S$ is. In particular it is a reduced ring. Fix $x\in X$ and a character $\chi$. If $S^{\chi}_{x}=0$, then, for every $1\leq b ,c \leq p-1$ such that $b c \equiv 1 \bmod p$ we have $(S^{\chi^b}_{x})^{c} \subseteq S^{\chi^{b c}}_{x}=0$. Since $S_{x}$ is reduced and $G$ is $p$-cyclic, it follows that $S^{\chi}_{x}=0$ for all nontrivial characters $\chi$. Thus $S_x=S^{\chi_0}_{x}=K_{x}$, contradicting the hypothesis that $S_{x}$ is $G$-Galois, which implies that $\dim_{K_{x}} S_{x} = \# G = p \geq 2$. Hence, we may assume that $S^{\chi}_{x}\neq 0$ for all nontrivial characters and choose $s\in S^{\chi}_{x}\setminus\{0\}$. Then, the map
\[
	S^{\chi}_{x} \xrightarrow{\cdot s^{p-1}}  (S^{\chi}_{x})^{p} \subset S^{\chi_0}_{x} = K_{x}
\]
sends $s$ to $s^p\in K_{x}$, which is invertible since $s \neq 0$ and $S_{x}$ is reduced. Hence $\dim_{K_{x}} S^{\chi}_{x}=1$.

\item $G$ acts faithfully since $S^{\chi}$ has rank $1$ and $\chi(g) - 1 \in \k^{*}$ for nontrivial $\chi,g$.

\item It remains to see that $g\in G\setminus\{\id\}$ and $\id$ are strongly distinct. Since $S_{x}$ is a finite reduced $K_{x}$-algebra, it is a finite product of finite field extensions of $K_{x}$:
\begin{equation}
\label{E:Sx field decomposition}
	S_{x}  = \prod_{j} K_{x,j}.
\end{equation}
This is an isomorphism of $K_{x}$-algebras mapping $K_{x}$ to the diagonal on the right hand side (analogously to the situation we saw in Lemma~\ref{L:structureKxntx} of Section~\ref{subsec:isomorphisms of adelic algebras}).

Recalling that $S$ injects into $\prod_{x} S_{x}$ (Lemma~\ref{L:S injects prod Sx}), we see that to prove that the elements of $G$ are strongly distinct, it suffices to consider the idempotents $1_{x,j} \in S$ corresponding to a $1$ in the factor $K_{x,j}$ of $S_{x}$ and $0$ elsewhere.

Suppose on the contrary that $g$ and $\id$ are not strongly distinct. Then for some idempotent of the form $1_{x,j}$ we have $g(s) 1_{x,j} = s 1_{x,j}$ for all $s \in S$. Since $G$ is $p$-cyclic, this relation in fact holds for all $g \in G$. This implies that the factor $K_{x,j} \subseteq S_{x}^{G} = K_{x}$ via~\eqref{E:Sx field decomposition}, and this can only happen if it is the only factor, that is to say, $S_{x} = K_{x}$, which since $\# G \geq 2$, contradicts that $S_{x}$ is $G$-Galois over $K_{x}$.
\qedhere
\end{itemize}
\end{proof}

\begin{remark}
\label{R:Galois rank general}
Theorem~\ref{T:galois local global} thus establishes a local-global criterion for Galois extensions of $\A_{X}$. Along with Theorem~\ref{T:Pic adeles} proving the triviality of the Picard group of $\A_{X}$, we obtain the freeness of the $\chi$-eigenspaces. 

In general, for a $p$-cyclic $G$-Galois extension of rings $S/R$ where $R$ is $p$-Kummerian, and a character $\chi : G \to R^{*}$, it is known that the summands $S^{\chi}$ are locally free $R$-modules of rank $1$, they satisfy $S^{\chi} S^{\psi} = S^{\chi \psi}$, and $(S^{\chi})^{p}$ is free. If in addition the $p$-torsion part of the Picard group of the ring $R$ is trivial, i.e. $\Pic(R)[p] = 0$, then the $S^{\chi}$ are free (\cite[\S 11, 13]{Borevich} or \cite[Propositions 5.3, 5.4, 6.5, Ch. 0]{Greither})).
\end{remark}

% (end)
% end subsection structure of p cyclic extensions

%%%%%%%%%%%%%%%%%%%%%%%%%%%%%%%%%%%%%%%%%%%%%%%%%%%%%%%%%%%%%%%%%%%%%%%%%%%%%%%%%%%
%%%%%%%%%%%%%%%%%%%%%%%%%%%%%%%%%%%%%%%%%%%%%%%%%%%%%%%%%%%%%%%%%%%%%%%%%%%%%%%%%%%
\subsection{\texorpdfstring{Classification of $\cc_{p}$-Galois extensions}{Classification of Cp-Galois extensions}}
\label{subsec:classification of cc_ p galois extensions}
% (fold)
%%%%%%%%%%%%%%%%%%%%%%%%%%%%%%%%%%%%%%%%%%%%%%%%%%%%%%%%%%%%%%%%%%%%%%%%%%%%%%%%%%%
%%%%%%%%%%%%%%%%%%%%%%%%%%%%%%%%%%%%%%%%%%%%%%%%%%%%%%%%%%%%%%%%%%%%%%%%%%%%%%%%%%%

We begin by adapting the notion of primitive element to the special case of Galois ring extensions of $\A_{X}$. The first important result is Theorem~\ref{T:primitive element theorem p-cyclic}, which may be viewed as a version of Hilbert's Theorem 90 for $p$-cyclic ring extensions of so-called Kummerian rings.

It should be noted that, in our framework, the correct definition and the existence of primitive elements is not clear. Although there are some generalizations to certain kinds of commutative rings of the classical statement for field extensions, to our knowledge, none of these can be applied in our setup. For instance, the results of~\cite{AndersonDobbsMullins} can not be applied since $\A_{X}$ has infinitely many idempotents.

An essential tool in the study of $p$-cyclic Galois ring extensions of $\A_{X}$ is the Harrison group $\H(\A_{X},\cc_{p})$, whose general definition we recall below (Definition~\ref{D:Harrison set (R,G)}). We determine the structure of this group, thus completely classifying the $p$-cyclic Galois extensions of $\A_{X}$, in Theorem~\ref{T:classification Cp-extensions AX general}.

We bring to mind some further key concepts from the Galois theory of rings as presented in~\cite{Borevich,Harrison,Greither}. Firstly, recall the notion of equivariance, which makes sense in our context of ring extensions $S/R$ under the action of a group by $R$-automorphisms.

\begin{definition}[Equivariant Isomorphism]
\label{D:equivariantly isomorphic extensions}
Let $R$ be a commutative ring and $G$ a fixed finite group. For $i=1,2$, let $S_{i}$ be a ring extension of $R$ with a faithful action of $G$ by $R$-automorphisms of $S_{i}$. We say that the pairs $(S_{1},G)$ and  $(S_{2},G)$ are equivariantly isomorphic, or simply $G$-isomorphic, via $\phi$, if $\phi$ is an $R$-algebra isomorphism $\phi: S_{1} \isomto S_{2}$ such that $\phi \circ g = g \circ \phi$ for all $g\in G$. 
\end{definition}

It is clear that a $G$-isomorphism preserves the $G$-Galois property of a ring extension. Harrison showed how to classify the set of  Galois extensions of a given ring $R$ and group $G$ with a fixed action, via the following object.

\begin{definition}[The Harrison group~\cite{Harrison}]
\label{D:Harrison set (R,G)}
Given a group $G$, the set of all $G$-isomorphism classes of $G$-Galois ring extensions $S$ over a fixed base ring $R$ with a fixed faithful action of $G$ is called the \emph{Harrison set} of $(R,G)$ and denoted by $\Har(R,G)$. When $G$ is a finite abelian group, $\Har(R,G)$ can be endowed with a group structure. In this case it is called the \emph{Harrison group}.
\end{definition}

The general definition of the group law is summarized in~\cite[p.10]{Greither}. With regard to the decomposition~\eqref{E:eigenspace decomposition}, the product of two $G$-Galois extensions $S_{i}/R$ for $i=1,2$ is given by
\begin{equation}
\label{E:Harrison product wrt eigenspaces}
	S_{1} \cdot S_{2} \eqdef \bigoplus_{\chi} (S_{1}^{\chi} \otimes S_{2}^{\chi}),
\end{equation}
where $G$ acts on the summand $S_{1}^{\chi} \otimes S_{2}^{\chi}$ via $g(s_{1} \otimes s_{2}) \eqdef \chi(g)(s_{1} \otimes s_{2})$. One checks that this product factors through equivariant equivalence and thus defines the group law on the Harrison group $\H(R,G)$.

The neutral element with respect to this product is the so-called trivial $G$-Galois extension, defined by $R^{(G)} \eqdef \bigoplus_{\chi} R$, considered as the set of maps from $G$ to $R$ under the standard sum and product, and with the action of $G$ given by $g((r_{\chi})_{\chi}) = (r_{g^{*}\chi})_{\chi}$ where $g^{*}$ denotes composition with multiplication by $g$.

With a view towards determining the structure of the Harrison group $\H(\A_{X},\cc_{p})$, we will give a series of auxiliary results.

Consider the case when $S/R = L/K$ is a $p$-cyclic Galois extension of fields, where $\chr(K) \neq p$ and $K$ contains the $p$th roots of unity (the $p$-Kummerian condition). Let $G = \Gal(L/K)$ and $\chi : G \to \mu_{p}(K^{*})$ a nontrivial character. If $g$ is a generator of $G$ then $\chi(g)$ is a primitive $p$-th root of unity. Since $\mathbf{N}_{L/K}(\chi(g)) = \chi^{p}(g) = 1$, Hilbert's Theorem 90 implies that there is an element $\alpha \in L^{*}$ with $\alpha \notin K$, such that $\chi(g) = g(\alpha)/\alpha$. In other words,  there is a nontrivial element $\alpha \in L^{\chi}$. Now $u \eqdef \alpha^{p} \in K^{*}$ and thus $L = K(\alpha) = K[T]/(T^{p} - u)$ is a Kummer extension and $\alpha$ is a primitive element in the classical sense of field theory.

This situation leads us to the following definition in the general context of ring extensions.

\begin{definition}[$G$-primitive element]
\label{D:G-primitive}
Let $R$ be a $p$-Kummerian ring and $S$ an $R$-algebra on which a $p$-cyclic group $G$ acts via $R$-automorphisms. Given a \emph{nontrivial} character $\chi : G \to \mu_{p} \subseteq R^{*}$, an element $\alpha \in S$ is called $(G,\chi)$-primitive if:
\begin{itemize}

\item $1, \alpha, \dots, \alpha^{p-1}$ is an $R$-module basis of $S$.

\item $g(\alpha) = \chi(g) \alpha$ for $g \in G$, i.e. $\alpha \in S^{\chi}$.

\end{itemize}
In this case we also say that $\alpha$ is $G$-primitive with character $\chi$.
% In view of Remark~\ref{R:G-primitive element fields}, we may also use the same terminology in the case of a $p$-cyclic field extension $L/K$.
\end{definition}

Note that the existence of a $G$-primitive element implies that the group acts faithfully, because if $\zeta \in \mu_{p} \setminus \{1\}$ then $\zeta - 1 \in R^{*}$.
Thus $g(s) = s$ for all $s \in S$ iff $g(\alpha) = \alpha$ iff $(\chi(g) - 1)\alpha = 0$, and $\chi(g) \in \mu_{p} \setminus \{1\}$ if $g \neq 1$.

In what follows, since we consider free $R$-algebras $S$ of rank $n$, the notion of characteristic polynomial is defined for elements $a \in S$ as usual, namely, $C_{a}(T) \in R[T]$ is the determinant of the map given by multiplication by $T - a$ on $S[T]$. Note that minimal polynomials need not exist for algebraic elements, when dealing with rings with zero divisors (such as $\A_{X}$).

% In particular, the standard notion of degree as the largest power appearing with a nonzero coefficient, one has $\deg(PQ) \leq \deg(P) + \deg(Q)$ and strict inequality may hold. As a consequence, minimal polynomials may not exist for algebraic elements.
% To avoid confusion, we will only use this terminology for cases when the base ring is a field.

\begin{proposition}
\label{P:G-primitive equivalences}
Let $S$ be a separable algebra over a $p$-Kummerian ring $R$ on which the $p$-cyclic group $G$ acts via $R$-automorphisms with $S^{G} = R$. Fix a nontrivial character $\chi : G \to \mu_{p} \subseteq R^{*}$. 
For an element $\alpha \in S^{\chi}$, the following are equivalent:
\begin{enumerate}

\item $\alpha$ is $(G,\chi)$-primitive, i.e. $1, \alpha, \alpha^{2}, \dots, \alpha^{p-1}$ is an $R$-module basis of $S$.

% \item $\alpha^{b}$ generates $S^{\chi^{b}}$ as an $R$-module for all $b = 0,1,\dots, p-1$.

\item $\alpha^{p} \in R^{*}$.

\item $\alpha$ is invertible in $S$.

\end{enumerate}
If this is the case, then:
\begin{enumerate}

\setcounter{enumi}{3}

\item $\alpha$ generates $S^{\chi}$ as an $R$-module, and $S$ is a free $R$-module of rank $p$.

\item The characteristic polynomial of $\alpha$ is $C_{\alpha}(T) = T^{p} -\alpha^{p} \in R[T]$. It is separable and generates $\Ann(\alpha)$.

\item $S$ is equivariantly isomorphic to the $(G,\chi)$-Kummer extension (Definition~\ref{D:G-chi-Kummer extension}) $(\gkummer[p]{R}{u},G,\chi)$ for $u = \alpha^{p}$.
\end{enumerate}
\end{proposition}

\begin{proof}
\leavevmode	
\begin{itemize}

\item $(1) \implies (2)$: By (1) $S$ is a free $R$-module of rank $p$, so that it makes sense to talk about the characteristic polynomial. Noting that $\alpha^{p} \in S^{\chi^{p}} = S^{\chi_{0}} = S^{G} = R$, a straightforward computation shows that the matrix of multiplication by $\alpha$ with respect to the basis $\{\alpha^{b}\}_{b=0}^{p-1}$ is $\begin{psmallmatrix} 0 & \alpha^{p} \\ I_{p-1} & 0 \end{psmallmatrix}$ and that $C_{\alpha}(T) \eqdef \det(T - \alpha) = T^{p} - \alpha^{p} \in R[T]$. No monic polynomial over $R$ of degree less than $p$ annihilates $\alpha$, hence $\Ann(\alpha)=(C_{\alpha}(T))$ and $S \simeq R[T]/(C_{\alpha}(T))$. Since $S$ is separable, $C_{\alpha}(T)$ is a separable polynomial, and hence, by the same reasoning as in Proposition~\ref{P:Axt-separable}, replacing $\A_{X}$ with $R$ and the parameter vector $\t \in \A_{X}$ by $\alpha^{p} \in R$, we conclude that in fact $\alpha^{p} \in R^{*}$.

\item $(2) \implies (3)$: Trivial.

\item $(3) \implies (1)$: If $\alpha \in S^{\chi}$ and $\alpha \in S^{*}$ (which is certainly the case if $\alpha^{p} \in R^{*}$), then $\alpha^{b}$ generates $S^{\chi^{b}}$ for any $b$, since for any $\beta \in S^{\chi^{b}}$, we may write $\beta = (\beta/\alpha^{b}) \alpha^{b}$. Now $g(\beta/\alpha^{b}) = (\chi^{b}(g) \beta)/(\chi^{b}(g) \alpha^{b}) = \beta/\alpha^{b}$, therefore $\beta/\alpha^{b} \in S^{G} = R$. Hence $S = \bigoplus_{\chi} S^{\chi} = \bigoplus_{b=0}^{p-1} R \alpha^{b}$.
\end{itemize}
We have already seen that (4) and (5) hold. By (5) we conclude that $\gkummer[p]{R}{u} \eqdef R[T]/(T^{p} - u)$ where $u \eqdef \alpha^{p} \in R^{*}$. Thus there is a unique $R$-algebra morphism $\phi : \gkummer[p]{R}{u} \to S$ sending $T$ to $\alpha$, which is in fact an equivariant isomorphism between the  $(G,\chi)$-Kummer extension $(\gkummer[p]{R}{u},G,\chi)$ and $S$.
\end{proof}

The condition on primitive elements given in Proposition~\ref{P:G-primitive equivalences}(5) is consistent with~\cite[Theorem 3.4]{Nagahara} which relates the existence of a primitive element in a ring extension to the existence of an embedding into a Galois extension of the base ring.

The invertibility condition in Proposition~\ref{P:G-primitive equivalences}(3) for
$G$-primitive elements in a $p$-cyclic Galois ring extension of $\A_{X}$ with group $G$ is easily seen to be related to the criterion given in~\cite[Corollary 2.2]{Paques} for an element in an arbitrary Galois extension of rings to be primitive.

As we now show, the obstruction to the existence of a $G$-primitive element is essentially that the extension be $G$-Galois, under a certain additional condition on the ring $R$. The relations between the various hypotheses is made clear below in Remark~\ref{R:Kummer sequence general}.

\begin{theorem}[Primitive Element Theorem for $p$-cyclic ring extensions]
\label{T:primitive element theorem p-cyclic}
Let $R$ be a $p$-Kummerian ring with $\Pic(R)[p] = 0$ and $S$ a separable $R$-algebra on which a $p$-cyclic group $G$ acts via $R$-automorphisms, with $S^{G} = R$. Fix a nontrivial character $\chi : G \to \mu_{p} \subseteq R^{*}$. Then the following are equivalent:
\begin{enumerate}

\item $(S,G)$ is a $G$-Galois extension of $R$.

\item $S$ has a $(G,\chi)$-primitive element.

\end{enumerate}
If this is the case,
\begin{enumerate}

\setcounter{enumi}{2}

\item $S^{\chi}$ is free of rank $1$ over $R$. Any generator $\alpha$ of $S^{\chi}$ is $(G,\chi)$-primitive.

\item $S$ is $G$-isomorphic to the $(G,\chi)$-Kummer extension $(\gkummer[p]{R}{u},G,\chi)$ with $u = \alpha^{p}$, for any $(G,\chi)$-primitive element $\alpha$.

\item The quotient of two $(G,\chi)$-primitive elements lies in $R^{*}$.

\item If $\alpha$ is $(G,\chi)$-primitive then for $b \in \Zps$, the $(G,\chi^{b})$-primitive elements are those of the form $u \alpha^{b}$ where $u \in R^{*}$. 

\item Given nontrivial characters $\chi_{i}$ for $i=1,2$, with $\chi_{1} = \chi_{2}^{b}$ for some $b \in \Zps$, any $(G,\chi_{i})$-primitive elements $\alpha_{i}$ satisfy
$\alpha_{1}/\alpha_{2}^{b} \in R^{*}$.

\end{enumerate}
\end{theorem}

\begin{proof}
If $S/R$ is $G$-Galois, then $S^{\chi}$ is a free $R$-module of rank $1$ and $S^{\chi} S^{\psi} = S^{\chi \psi}$ (see Remark~\ref{R:Galois rank general}). Choose any generator $\alpha$ of $S^{\chi}$. For $0 \leq b \leq p-1$ we have $S^{\chi^{b}} = (S^{\chi})^{b} = R \alpha^{b}$, hence $S = \bigoplus_{b=0}^{p-1} S^{\chi^{b}} = \bigoplus_{b=0}^{p-1} R \alpha^{b}$. Thus (2), (3) and (4) follow from
Proposition~\ref{P:G-primitive equivalences}.

Conversely, if $\alpha$ is $(G,\chi)$-primitive, then by Proposition~\ref{P:G-primitive equivalences} we know that $u \eqdef \alpha^{p} \in R^{*}$ and $S$ is $G$-isomorphic to $(\gkummer[p]{R}{u},G,\chi)$, and in particular $S$ is a $G$-Galois extension of $R$.

Since the nontrivial character $\chi$ was arbitrary, Proposition~\ref{P:G-primitive equivalences} and (3) say that for $1 \leq b \leq p-1$, the $(G,\chi^{b})$-primitive elements are exactly the generators of $S^{\chi^{b}}$. Thus (5)--(7) follow straightforwardly.
\end{proof}

\begin{proposition}
\label{P:isomorphism criterion p-cyclic}
Let $R$ be a $p$-Kummerian ring, $G$ a $p$-cyclic group, and $S_{1}$, $S_{2}$ two $G$-Galois extensions of $R$. Fix a nontrivial character $\chi : G \to \mu_{p} \subseteq R^{*}$ and $(G,\chi)$-primitive elements $\alpha_{i} \in S_{i}^{\chi}$. Then the following are equivalent:
\begin{enumerate}

\item $S_{1}$ and $S_{2}$ are $G$-equivariantly isomorphic.

\item $\alpha_{1}^{p}/\alpha_{2}^{p} \in R^{*p}$.

\end{enumerate}
\end{proposition}

\begin{proof}
Suppose that $\phi : S_{1} \to S_{2}$ is a $G$-equivariant $R$-isomorphism. Then $\phi(\alpha_{1}) \in S_{2}^{\chi} = u \alpha_{2}$ for some $u \in R^{*}$. Raising to the $p$th power and keeping in mind that $\alpha_{i}^{p} \in R^{*}$ we obtain $\alpha_{1}^{p} = u^{p} \alpha_{2}^{p}$.

Conversely, suppose that (2) holds with $\alpha_{1}^{p} = u^{p} \alpha_{2}^{p}$ for some $u \in R^{*}$. Let $u_{i} \eqdef \alpha_{i}^{p}$. By Theorem~\ref{T:primitive element theorem p-cyclic} we know that $S_{i}$ is $G$-equivariantly isomorphic to the $(G,\chi)$-Kummer extension $\gkummer[p]{R}{u_{i}} = R[T]/(T^{p} - u_{i})$. It is now straightforward to check that $T \mapsto u T$ defines a $G$-equivariant $R$-isomorphism between $\gkummer[p]{R}{u_{1}}$ and $\gkummer[p]{R}{u_{2}}$.
\end{proof}

\begin{proposition}
\label{P:Kummer sequence left}
For a $p$-Kummerian ring $R$, a $p$-cyclic group $G$, and a choice of nontrivial character $\chi$, the map $i = i_{(R,\chi)}$ sending $u \in R^{*}$ to the $(G,\chi)$-Kummer extension $\gkummer[p]{R}{u}$ induces an injective group homomorphism
\[
	\faktor{R^{*}}{R^{*p}} \hookrightarrow \H(R,G),
\]
which is an isomorphism if $\Pic(R)[p] = 0$. We shall refer to this map as the \emph{Kummer map} corresponding to $(R, \chi)$.
\end{proposition}

\begin{proof}
We check that $i$ is a group homomorphism, i.e. that the Harrison product of two $(G, \chi)$-Kummer extensions $R{u_i}$ for $i=1,2$ is $R{(u_1 u_2)}$. Fixing $(G,\chi)$-primitive elements $\alpha_i\in R{u_i}^{\chi}$ (which exist by Theorem~\ref{T:primitive element theorem p-cyclic}) it is clear that $\alpha_{1} \otimes \alpha_{2}$ belongs to $R{u_1}^{\chi} \otimes R{u_2}^{\chi}$, which by~\eqref{E:Harrison product wrt eigenspaces} is a direct summand of $ R{u_1} \cdot R{u_2}$, and $(\alpha_{1} \otimes \alpha_{2})^{p} \in R^{*}$. Hence, by Proposition~\ref{P:G-primitive equivalences}, $R{u_1} \cdot R{u_2}$ is equivariantly isomorphic to $R{(u_1 u_2)}$. By Proposition~\ref{P:isomorphism criterion p-cyclic}, $\ker(i) = R^{*p}$.

When $\Pic(R)[p] = 0$, the surjectivity of the Kummer map is given in Theorem~\ref{T:primitive element theorem p-cyclic}(4).
\end{proof}

\begin{remark}[The general Kummer sequence]
\label{R:Kummer sequence general}
Proposition~\ref{P:Kummer sequence left} is a special case of the following general fact. For a $p$-Kummerian ring, $R$, the sequence of groups:
\begin{equation}
\label{E:Kummer sequence general}
	1 \xrightarrow{} \faktor{R^{*}}{R^{*p}}
	  \xrightarrow{i_{(R,\chi)}} \H(R,\cc_{p})
	  \xrightarrow{} \Pic(R)[p]
	  \xrightarrow{} 1
\end{equation}
is exact. This sequence is called the \emph{Kummer sequence}.

Thus if $\Pic(R)$ has no $p$-torsion, the Harrison group $\H(R,G)$ is isomorphic to $R^{*}/R^{*p}$.
% , which opens up the possibility of applying or adapting some ideas from class field theory and reciprocity laws to this context.

As regards Theorem~\ref{T:primitive element theorem p-cyclic}, in general, without assuming $\Pic(R)[p]= 0$, the existence of a $G$-primitive element is equivalent to the existence of a \emph{normal basis} (see e.g.~\cite[\S13]{Borevich} for the  definition in the context of ring extensions). Indeed, one may show that under the conditions of Proposition~\ref{P:G-primitive equivalences}, an element $\alpha \in S^{\chi}$ is $(G,\chi)$-primitive if and only if $1 + \alpha + \alpha^{2} + \cdots + \alpha^{p-1}$ generates a normal basis (this is essentially~\cite[Proposition 6.5]{Greither}).
\end{remark}

In the case of the adele ring $R = \A_{X}$, we can easily understand the quotient $R^{*}/R^{*p}$.

\begin{definition}
\label{D:idele valuation vector}
For an idele $u \in \I_{X}$, define its valuation vector as the element
\begin{equation}
\label{E:idele valuation vector}
	\upsilon(u) \eqdef (\upsilon_{x}(u_{x}) \bmod p)_{x} \in  \bigoplus_{x \in X} \Zp.
\end{equation}
\end{definition}

Note that~\eqref{E:idele valuation vector} is a group homomorphism and by~\eqref{E:n th powers valuation 0 mod n} we have
\begin{equation}
\label{E:IXp kernel upsilon}
	\I_{X}^{p} = \ker \upsilon : \I_{X} \to \bigoplus_{x} \Zp.
\end{equation}

\begin{theorem}[Classification of $p$-cyclic Galois extensions of $\A_{X}$]
\label{T:classification Cp-extensions AX general}
Let $X$ be a projective, irreducible, non-singular curve over an algebraically closed field $\k$ and $\Sigma$ the function field of $X$. Fix a nontrivial character $\chi : \cc_{p} \to \k^{*}$, where $p$ is a prime different from $\chr \Bbbk$. Then the following diagram is commutative
\begin{equation}
\label{E:Kummer sequence local global}
\begin{tikzcd}
	    \bigoplus\limits_{x \in X} \faktor{\Z}{(p)}
		\arrow[d, hook] 
	  & \faktor{\I_{X}}{\I_{X}^{p}}
	    \arrow[r, "i_{(\A_{X},\chi)}", "\sim"'] \arrow[l, "\upsilon"',"\sim"]  \arrow[d, hook] 
	  & \H(\A_{X},\cc_{p}) 
	    \arrow[d, hook]
	\\
	    \prod\limits_{x \in X} \faktor{\Z}{(p)}
	  & \prod\limits_{x \in X} \faktor{K_{x}^{*}}{K_{x}^{*p}}
	    \arrow[r, "(i_{(K_{x},\chi)})_{x}", "\sim"'] \arrow[l, "(\upsilon_{x})_{x}"' ,"\sim"]
	  & \prod\limits_{x \in X} \H(K_{x},\cc_{p})
\end{tikzcd}
\end{equation}
where $u\in \I_{X}$ (respectively, $u \in K_{x}^{*}$) is mapped to the class of the $(\cc_{p}, \chi)$-Kummer extension $(\AXn[p]{u},\cc_{p},\chi)$ (respectively, $(\gkummer[p]{K_{x}}{u},\cc_{p},\chi)$).
\end{theorem}

\begin{proof}
First, we deal with the top row. As we saw in Proposition~\ref{P:Kummer sequence left}, the 
map $\I_X \to  \H(\A_{X},\cc_{p})$ is an injective group homomorphism, and by Theorem~\ref{T:Pic adeles} and~\ref{T:primitive element theorem p-cyclic}, the map is surjective, thus inducing the top right isomorphism.

Similar reasoning shows that there are group isomorphisms $K_{x}^{*}/K_{x}^{*p} \simeq \H(K_{x},\cc_{p})$ for all $x\in X$. 

It is clear that the left horizontal arrows, induced by the valuations, are isomorphisms, and that the first two vertical arrows are injective. The existence of the right vertical arrow is a consequence of Theorem~\ref{T:galois local global}, and it is straightforward to check commutativity. Finally, the injectivity on the right follows from the rest.
\end{proof}

\begin{definition}[Valuation vector of a $p$-cyclic Galois extension]
\label{D:valuation vector S G chi}
Let $(S,G)$ be a $p$-cyclic Galois ring extension of $\A_{X}$ and $\chi : G \to \mu_{p} \subseteq \k^{*}$ a nontrivial character, the \emph{valuation vector} associated to the triple $(S,G,\chi)$ is the image of the extension under the map $\upsilon \circ i_{(\A_{X},\chi)}^{-1}$. We will denote it by $\upsilon(S,G,\chi)$.  Chasing the various definitions, one easily sees that
\begin{equation}
\label{E:(G,chi)-valuation vector}
	\upsilon(S,G,\chi) = \upsilon(\alpha^{p}) = (\upsilon_{x}(\alpha^{p}) \bmod p)_{x} \in \bigoplus_{x \in X} \Zp,
\end{equation}
where $\alpha$ is a $(G,\chi)$-primitive element, which exists by Theorem~\ref{T:primitive element theorem p-cyclic}. This same theorem shows that the above expression does not depend on the choice of $\alpha$.
\end{definition}

The top row of~\eqref{E:Kummer sequence local global} shows that $\upsilon(S,G,\chi)$ is invariant under $G$-isomorphism.

It is clear that for a $(G,\chi)$-Kummer extension (Definition~\ref{D:G-chi-Kummer extension}) $(\AXn[p]{\t},G,\chi)$ of $\A_{X}$, we have
\begin{equation}
\label{E:valuation vector (G,chi)-Kummer}
		\upsilon(\AXn[p]{\t},G,\chi) = \upsilon(\t),
\end{equation}
since the condition $g(T) = \chi(g) T$ means that the class of $T$ is a $(G,\chi)$-primitive element $\alpha$ satisfying $\alpha^{p} = \t$.

Now, consider any $p$-cyclic Galois extension $(\AXn[p]{\t},G)$ of $\A_{X}$ and a nontrivial character $\chi : G \to \mu_{p}(\k^{*})$, not necessarily a $(G,\chi)$-Kummer extension, namely the class of $T$, although it is a primitive element whose $p$th power is an idele, need not satisfy $g(T) = \chi(g) T$.
Nevertheless, there are $(G,\chi)$-primitive elements. If $\alpha$ is one such, by Theorem~\ref{T:primitive element theorem p-cyclic}(4), $(\AXn[p]{\t},G,\chi)$ is $G$-isomorphic to the $(G,\chi)$-Kummer extension $(\AXn[p]{(\alpha^{p})},G,\chi)$. By~\eqref{E:valuation vector (G,chi)-Kummer} its valuation vector is $\upsilon(\alpha^{p})$. In particular, since $\AXn[p]{\t}$ and $\AXn[p]{(\alpha^{p})}$ are isomorphic $\A_{X}$-algebras, by Theorem~\ref{T:AXt isomorphic iff ram equal and idele ratio} we conclude that 
\begin{equation}
\label{E:ram ap = ram t}
	\Ram(\alpha^{p}) = \Ram(\t),
\end{equation}
that is to say, $\upsilon_{x}(\alpha^{p}) \equiv 0 \bmod p$ if and only if $x \notin \Ram(\t)$.
In particular, $\upsilon(\AXn[p]{\t},G,\chi)$ is the null vector if and only if $\Ram(\t) = \emptyset$.

% \begin{remark}[Induced topology]
% \label{R:topology on p-cyclic S}
% Combining Theorems~\ref{T:primitive element theorem p-cyclic} and~\ref{T:classification Cp-extensions AX general} with Corollary~\ref{C:algebraic isomorphism AXT is topological}, we conclude that any $\cc_{p}$-Galois ring extension $S$ of $\A_{X}$ has a natural topology, induced by choosing any $\cc_{p}$-isomorphism of $S$ with a $(\cc_{p},\chi)$-Kummer adelic algebra $(\AXn[p]{\t},\cc_{p},\chi)$.
% \end{remark}

Theorem~\ref{T:classification Cp-extensions AX general} says that the $\cc_{p}$-Galois extensions of the adele ring $\A_{X}$ are, up to equivariant isomorphism, the $(\cc_{p},\chi)$-Kummer adelic extensions $(\AXn[p]{\t},\cc_{p},\chi)$. This may be regarded as an analog for the base ring $R = \A_{X}$ of the classical result for fields stating that a $p$-cyclic Galois field extension $L/K$, where $\mu_{p} \subseteq K^{*}$, is a radical extension of the form $L = K(u^{1/p})$ where $u \in K^{*} \setminus K^{*p}$ (see~\cite[Ch. III, \S2]{CaFr}). The dependence of the isomorphism on the choice of $\cc_{p}$-primitive element was explored in Proposition~\ref{P:isomorphism criterion p-cyclic}.

% (end)
% subsection classification of cc_ p galois extensions end

%%%%%%%%%%%%%%%%%%%%%%%%%%%%%%%%%%%%%%%%%%%%%%%%%%%%%%%%%%%%%%%%%%%%%%%%%%%%%%%%%%%
%%%%%%%%%%%%%%%%%%%%%%%%%%%%%%%%%%%%%%%%%%%%%%%%%%%%%%%%%%%%%%%%%%%%%%%%%%%%%%%%%%%
\subsection{\texorpdfstring{Conjugacy of $p$-cyclic Galois extensions}{Conjugacy of p-cyclic Galois extensions}}
\label{subsec:conjugacy of p cyclic galois extensions}
% (fold)
%%%%%%%%%%%%%%%%%%%%%%%%%%%%%%%%%%%%%%%%%%%%%%%%%%%%%%%%%%%%%%%%%%%%%%%%%%%%%%%%%%%
%%%%%%%%%%%%%%%%%%%%%%%%%%%%%%%%%%%%%%%%%%%%%%%%%%%%%%%%%%%%%%%%%%%%%%%%%%%%%%%%%%%

As we have mentioned before, in the Galois theory of commutative rings, there can be many ``Galois groups'', as opposed to what happens for field extensions. 

In terms of ring extensions, this entails considering different, although isomorphic, Galois groups of a fixed extension. This is essentially the same as considering a fixed group $G$ and varying its action, leading to a notion of conjugacy for ring extensions (Definition~\ref{D:conjugate extensions}) generalizing that of equivariant isomorphism.

Later on in \S\ref{subsec:galois equivalent subgroups of automorphisms} we will turn our attention to the relation between conjugate extensions and group-theoretical properties.

% We now introduce the concept of conjugacy for $G$-ring extensions and
% % , which  generalizes the usual notion of $G$-isomorphism (Definition~\ref{D:equivariantly isomorphic extensions}).
% then give a classification, using certain equivalence classes of valuation vectors, of the $p$-cyclic Galois extensions of $\A_{X}$ up to conjugacy.

\begin{definition}[Conjugacy of $G$-ring extensions]
\label{D:conjugate extensions}
Let $R$ be a commutative ring. For $i=1,2$, let $S_{i}$ be a ring extension of $R$ with a faithful action of a group $G_{i}$ by $R$-automorphisms of $S_{i}$. We say that $(S_{1},G_{1})$ and  $(S_{2},G_{2})$ are conjugate via $(\phi,\tau)$ if $\phi: S_{1} \isomto S_{2}$ is an $R$-algebra isomorphism and $\tau : G_{1} \isomto G_{2}$ is a group isomorphism, such that $\phi\circ g = \tau(g) \circ \phi$ for all $g\in G_1$. We will denote this relation by $(S_{1}, G_{1}) \sim (S_{2}, G_{2})$.
\end{definition}

Note that the group isomorphism $\tau$ is in fact determined by the $R$-algebra isomorphism $\phi$, namely $\tau(g) = g^{\phi} \eqdef \phi \circ g \circ \phi^{-1}$, although it is convenient to denote it as part of a pair $(\phi,\tau)$ as we are doing here.

\begin{remark}
\label{R:conjugation preserves Galois}
If one looks at the definition, it is immediately clear that being a Galois extension $S/R$ is preserved by conjugation.
\end{remark}

\begin{definition}
\label{D:twist (S,G) aut}
Suppose $(S,G)$ is a ring extension of $R$ with $G$ a group acting faithfully by $R$-automorphisms on $S$. Given an automorphism $\tau \in \Aut(G)$ of $G$, we define the \emph{twist} of $(S,G)$ by $\tau$, denoted by $(S,G)^{\tau}$, to be the same ring extension $S/R$ but with the action of $G$ now given by $g(s) \eqdef \tau(g) s$.
\end{definition}

\begin{remark}
\label{R:Harrison bifunctoriality}
The above action on Galois extensions induced by an automorphism of the group is a special case of the more general result that $\H(R,\mbox{--})$ is a functor from finite abelian groups to abelian groups (\cite[Theorem 3.2]{Greither} or~\cite[\S1]{Harrison}). In particular, a group automorphism of $G$ induces a group automorphism of $\H(R,G)$.
\end{remark}

One immediately sees that $(S,G)^{\tau}$ is conjugate to $(S,G)$ via $(\id,\tau)$.
In particular, conjugacy classes of $p$-cyclic Galois ring extensions correspond bijectively to the quotient set of $\cc_{p}$-isomorphism classes of $\cc_{p}$-Galois extensions modulo $\Aut(\cc_{p})$, i.e.
\begin{equation}
\label{E:correspondence conjugacy p-cyclic mod aut}
	\left\{
	\begin{minipage}[c]{0.42\textwidth}
	\raggedright
	Conjugacy classes of $p$-cyclic Galois extensions $(S,G)$ of $R$
	\end{minipage}
	\right\}
	\correspondence
	\faktor{\H(R,\cc_{p})}{\Aut(\cc_{p})}.
\end{equation}
We may also consider the following action of $\Zps$ by automorphisms on the quotient $R^{*}/R^{*p}$ defined by
\[
	b \in \Zps \longmapsto (u \mapsto u^{b}) \in \Aut(R^{*}/R^{*p}).
\]
When $R$ is $p$-Kummerian, a choice of nontrivial character $\chi : \cc_{p} \to \mu_{p} \subseteq R^{*}$ gives us:
\begin{itemize}

\item The Kummer map $i_{(R,\chi)} : R^{*}/R^{*p} \to \H(R,\cc_{p})$ from Proposition~\ref{P:Kummer sequence left}.

\item The group isomorphism $\chi^{*} : \Aut(\mu_{p}) \simeq \Zps \isomto \Aut(\cc_{p})$ induced by $\chi : \cc_{p} \isomto \mu_{p} \subseteq R^{*}$.

\end{itemize}
The next result shows how all of these relate to one another.

\begin{lemma}
\label{L:Kummer map actions}
For a $p$-Kummerian ring $R$ and a choice of nontrivial character $\chi$ of $\cc_{p}$, we have:
\begin{enumerate}

\item $i_{(R,\chi)}(u^{b}) = i_{(R,\chi^{c})}(u)$, where $c$ is the inverse of $b$ modulo $p$.

\item $i_{(R,\chi \circ \tau)}(u) = (i_{(R,\chi)}(u))^{\tau^{-1}}$ for $\tau \in \Aut(\cc_{p})$.

\item $i_{(R,\chi)}$ is equivariant with respect to the actions of $\Zps$ and $\Aut(\cc_{p})$ via the group isomorphism $\chi^{*}$.

\end{enumerate}
\end{lemma}

\begin{proof}
(1) and (2) are straightforward computations using the definitions of the Kummer map and the twist by an automorphism. Using these relations, given $b \in \Zps$, let $\tau = \chi^{*}(b)$ and $bc \equiv 1 \bmod p$. Then $i_{(R,\chi)}(u^{b}) = i_{(R,\chi^{c})}(u) = i_{(R,\chi \circ \tau^{-1})}(u) = (i_{(R,\chi)}(u))^{\tau}$, which gives (3).
\end{proof}

\begin{theorem}
\label{T:correspondence conjugacy p-cyclic}
Let $R$ be a $p$-Kummerian ring with $\Pic(R)[p] = 0$. Then we have the following bijective correspondences:
\begin{equation}
\label{E:correspondence conjugacy p-cyclic triple}
	\begin{aligned}
	\faktor{(R^{*}/R^{*p})}{\Zps}
	&\correspondence
	\faktor{\H(R,\cc_{p})}{\Aut(\cc_{p})}
	\\
	&\correspondence
	\left\{
	\begin{minipage}[c]{0.42\textwidth}
	\raggedright
	Conjugacy classes of $p$-cyclic Galois extensions $(S,G)$ of $R$.
	\end{minipage}
	\right\}
	\end{aligned}
\end{equation}
In more explicit terms, given two $p$-cyclic Galois extensions $(S_{i}, G_{i})$ for $i = 1,2$, and choices $\alpha_{i}$ of $G_{i}$-primitive elements, the following are equivalent:
\begin{enumerate}
	
\item $(S,G_{1})$ and $(S,G_{2})$ are conjugate $p$-cyclic Galois extensions of $R$.

\item $\alpha_{1}^{p}/\alpha_{2}^{p b} \in R^{*p}$ for some $b \in \Zps$.

\end{enumerate}
\end{theorem}

\begin{proof}
By Proposition~\ref{P:Kummer sequence left}, the Kummer map $i_{(R,\chi)}$ is an isomorphism. Combining the correspondence~\eqref{E:correspondence conjugacy p-cyclic mod aut} and the equivariance of the $\Zps$ and $\Aut(\cc_{p})$-actions given in Lemma~\ref{L:Kummer map actions}, we conclude.

In terms of primitive elements, let $u_{i} \eqdef \alpha_{i}^{p}$, which as we saw in Proposition~\ref{P:G-primitive equivalences}, belongs to $R^{*}$. Then (2) says that $u_{1}$ and $u_{2}$ mod $R^{*p}$ lie in the same orbit under the action of $\Zps$. Hence the equivalence of (1) and (2) is immediate from~\eqref{E:correspondence conjugacy p-cyclic triple}.
\end{proof}

Note that in the commutative squares~\eqref{E:Kummer sequence local global}, the map $\upsilon$, which was defined in~\eqref{E:idele valuation vector} is also equivariant with respect to the action of $\Zps$ on $\bigoplus_{x} \Zp$ by multiplication. Thus \emph{every} map in the entire diagram is in fact equivariant.

This leads us to the following definitions.

\begin{definition}[Valuation class of a $p$-cyclic Galois extension]
\label{D:valuation class S G chi}
Let $(S,G)$ be a $p$-cyclic Galois ring extension of $\A_{X}$ and $\chi : G \to \mu_{p} \subseteq \k^{*}$ a nontrivial character, the \emph{valuation class} associated to the extension $(S,G)$ is the equivalence class of the valuation vector $\upsilon(S,G,\chi)$ under the action of $\Zps$:
\begin{equation}
\label{E:(S,G)-valuation class}
	\kappa(S,G) \eqdef [\upsilon(S,G,\chi)] \in 
	\biggl(\,\bigoplus_{x \in X} \faktor{\Z}{(p)}\biggr)
	\!\!\Bigm/\!\!
	\Bigl(\faktor{\Z}{(p)}\Bigr)^{*}.
\end{equation}
\end{definition}
It is clear that the valuation class $\kappa(S,G)$ does not depend on the choice of nontrivial character $\chi$.

For the base ring $R = \A_{X}$, Theorem~\ref{T:correspondence conjugacy p-cyclic} yields the following.

\begin{proposition}
\label{P:correspondence conjugacy p-cyclic AX}
The set of conjugacy classes of $p$-cyclic Galois extensions of the geometric adele ring $\A_{X}$ is in one-to-one correspondence with each of the following sets
\begin{equation}
\label{E:correspondence conjugacy p-cyclic AX}
	\biggl(\,\bigoplus_{x \in X} \faktor{\Z}{(p)}\biggr)
	\!\!\Bigm/\!\! \Bigl(\faktor{\Z}{(p)}\Bigr)^{*}
	\correspondence
	\faktor{(\I_{X}/\I_{X}^{p})}{\Zps}
	\correspondence
	\faktor{\H(\A_{X},\cc_{p})}{\Aut(\cc_{p})}
\end{equation}
Note that the composition of these correspondences from right to left is precisely the valuation class map $\kappa(S,\cc_{p})$.
\end{proposition}

Theorem~\ref{T:correspondence conjugacy p-cyclic} also holds when $R = K$ is a field of characteristic different from $p$ containing the $p$th roots of unity.

% (end)
% subsection conjugate subgroups end

%%%%%%%%%%%%%%%%%%%%%%%%%%%%%%%%%%%%%%%%%%%%%%%%%%%%%%%%%%%%%%%%%%%%%%%%%%%%%%%%%%%
%%%%%%%%%%%%%%%%%%%%%%%%%%%%%%%%%%%%%%%%%%%%%%%%%%%%%%%%%%%%%%%%%%%%%%%%%%%%%%%%%%%
\subsection{Local-global correspondence}
\label{subsec:local global correspondence}
% (fold)
%%%%%%%%%%%%%%%%%%%%%%%%%%%%%%%%%%%%%%%%%%%%%%%%%%%%%%%%%%%%%%%%%%%%%%%%%%%%%%%%%%%
%%%%%%%%%%%%%%%%%%%%%%%%%%%%%%%%%%%%%%%%%%%%%%%%%%%%%%%%%%%%%%%%%%%%%%%%%%%%%%%%%%%

In this section we consider the problem of relating the existence of a global conjugacy between two $p$-cyclic Galois extensions of the adele ring $\A_{X}$ to the existence of local conjugacies satisfying certain additional conditions. In other words, we seek an analog of Theorem~\ref{T:galois local global} for conjugacy rather than equivariant isomorphism.

We now introduce notation for both the full group of all $\A_{X}$-algebra automorphisms of an adelic algebra $\AXn[p]{\t}$,
\begin{equation}
\label{E:Gt}
	\G(\t) \eqdef \Aut_{\text{$\A_{X}$-alg}}\AXn[p]{\t}.
\end{equation}
%%
% Note that such automorphisms are automatically topological isomorphisms since we assume throughout that $\t \in \I_{X}$ (Corollary~\ref{C:algebraic isomorphism AXT is topological}).
as well as the corresponding group of $K_{x}$-algebra isomorphisms of $\Kxn[p]{t_{x}}$
at each point $x \in X$,
\begin{equation}
\label{E:Gtx}
	\G_{x}(\t) \eqdef \Aut_{\text{$K_{x}$-alg}}\Kxn[p]{t_{x}}.	
\end{equation}

\begin{proposition}
\label{P:structureGt}
For an idele $\t = (t_{x})$, we have
\begin{equation}
\label{E:Gtcomponents}
	\G(\t) = \prod_{x\in X} \G_{x}(\t)
\end{equation}
and
\begin{equation}
\label{E:Gtx structure}
		\G_{x}(\t)
		\simeq
		\begin{cases}
		\cc_{p} & \text{if $x\in\Ram(\t)$,} 
		\\
		\ss_p & \text{if $x\notin\Ram(\t)$}
	    \end{cases}
\end{equation}
for each $x \in X$, where $\cc_{n}$ and $\ss_{n}$ are respectively the cyclic and symmetric groups of degree $n$.
\end{proposition}

\begin{proof}
\eqref{E:Gtcomponents} follows immediately from Lemma~\ref{L:familiesOfIsomorphisms}.
For~\eqref{E:Gtx structure}, recall that by Lemma~\ref{L:structureKxntx}, when $x \notin \Ram(\t)$ we have $\Kxn[p]{t_{x}} \simeq \prod^{p} K_{x}$ and when $x \in \Ram(\t)$, $\Kxn[p]{t_{x}}/K_{x}$ is a $p$-cyclic field extension.
\end{proof}

For each field $K_{x}$, we have a pointwise analog of Proposition~\ref{P:correspondence conjugacy p-cyclic AX} which, since $\k$ is assumed algebraically closed of characteristic prime to $p$, may be seen to correspond to the formal analog of Puiseux's Theorem.

\begin{lemma}[Formal Puiseux's Theorem]
\label{L:Puiseux formal}
There is a unique abelian extension of $K_{x}$ of exponent $p$ in a fixed algebraic closure $K_{x}^{\text{a}}$, namely $E_{x} \eqdef K_{x}(K_{x}^{*(1/p)})$; it is a cyclic Kummer extension of degree $p$ generated by the $p$th root of any non $p$th power in $K_{x}$.
\end{lemma}

\begin{proof}
The hypotheses imply that $K_{x} \simeq \k\laurent{z_{x}}$ where $z_{x}$ is a uniformizer at $x$. Thus 
\begin{equation}
\label{E:Kummer Kx/Kxp=Z/p}
	  \ker(K_{x}^{*} \xrightarrow{\upsilon_{x}} \Z \xrightarrow{} \Zp)
	= \{ \lambda : \upsilon_{x}(\lambda) \equiv 0 \bmod p\} = K_{x}^{*p},
\end{equation}
therefore $K_{x}^{*}/K_{x}^{*p} \simeq \Zp$. The Kummer correspondence between subgroups of $K_{x}^{*}/K_{x}^{*p}$ and abelian extensions of $K_{x}$ of exponent $p$ yields existence and uniqueness. If $\Gamma$ is a subgroup of $K_{x}^{*} $ then either $\Gamma = K_{x}^{*p}$ or $K_{x}^{*p} \Gamma = K_{x}^{*}$, hence any non $p$th power generates the extension.
\end{proof}

The name of the theorem derives from the fact that in particular the unique extension may be realized as $\k\laurent{z_{x}^{1/p}}$ for a uniformizer $z_{x}$ (see~\cite[Ch. IV, Prop. 8]{SerreLocFields}).

The relation~\eqref{E:Kummer Kx/Kxp=Z/p} shows that for $K_{x}$, there are only two conjugacy classes, corresponding to the quotient $\Zp/\Zps$, namely:
\begin{itemize}

\item The class of the unique field extension given in the formal analog of Puiseux's Theorem.

\item The class of the neutral element in $\H(K_{x},\cc_{p})$, which is the $K_{x}$-algebra extension $\prod^{p} K_{x}$ where $\cc_{p}$ acts via a $p$-cycle.

\end{itemize}
Comparing this with Lemma~\ref{L:structureKxntx} which deals strictly with $K_{x}$-algebra isomorphisms, we obtain the following result.

\begin{proposition}
\label{P:local conjugacy Kxt}
For $t_{i,x} \in K_{x}^{*}$, $i=1,2$, let $(\Kxn[p]{t_{i,x}},G_{i,x})$
be $p$-cyclic Galois ring extensions of $K_{x}$. Then the following are equivalent:
\begin{enumerate}

\item $(\Kxn[p]{t_{1,x}},G_{1,x})$ and $(\Kxn[p]{t_{2,x}},G_{2,x})$ are conjugate Galois extensions.

\item $\Kxn[p]{t_{1,x}}$ and $\Kxn[p]{t_{2,x}}$ are isomorphic as $K_{x}$-algebras.

\end{enumerate}
\end{proposition}

\begin{proof}
\leavevmode
\begin{itemize}

\item $(1) \implies (2)$: This is trivial, by definition of conjugacy.

\item $(2) \implies (1)$: By Lemma~\ref{L:structureKxntx}, there are only two cases:

\end{itemize}

\begin{itemize}

\item Both $\Kxn[p]{t_{i,x}}$ are $p$-Kummer \emph{field} extensions of $K_{x}$. In this case, we conclude by Lemma~\ref{L:Puiseux formal}.

\item Both $\Kxn[p]{t_{i,x}}$ are isomorphic as $K_{x}$-algebras to $\prod^{p} K_{x}$. Fixing isomorphisms $\Kxn[p]{t_{i,x}} \isomto \prod^{p} K_{x}$, and generators $g_{i,x}$ of $G_{i,x}$, the latter correspond to $p$-cycles $\sigma_{i} \in \ss_{p}$. The isomorphism $\psi_{x} : \Kxn[p]{t_{1,x}} \to  \Kxn[p]{t_{2,x}}$ induced by $(a_{1}, \dots, a_{p}) = (a_{\sigma_{x}(1)}, \dots, a_{\sigma_{x}(p)})$, where $\sigma_{x} = \sigma_{2}^{-1} \sigma_{1}$ yields the desired conjugation.
\qedhere
\end{itemize}
\end{proof}

\begin{proposition}
\label{P:local conjugacy AXt}
For $\t_{i} \in \I_{X}$, $i=1,2$, let $(\AXn[p]{\t_{i}},G_{i})$ with $G_{i} \subseteq \G(\t_{i})$
be $p$-cyclic Galois ring extensions of $\A_{X}$. Then the following are equivalent:
\begin{enumerate}

\item $(\Kxn[p]{t_{1,x}},G_{1})$ and $(\Kxn[p]{t_{2,x}},G_{2})$ are conjugate Galois extensions for all $x \in X$, where $G_{i}$ acts on $\Kxn[p]{t_{i,x}}$ via projection.

\item $\AXn[p]{\t_{1}}$ and $\AXn[p]{\t_{2}}$ are isomorphic as $\A_{X}$-algebras.

\end{enumerate}
\end{proposition}

\begin{proof}
\leavevmode
\begin{itemize}

\item $(1) \implies (2)$: If $\Kxn[p]{t_{1,x}},G_{1})$ and $(\Kxn[p]{t_{2,x}},G_{2})$ are conjugate via $(\phi_{x}, \tau_{x} )$, then $\phi_{x}$ induces an isomorphism $\phi : \prod_{x\in X} \Kxn[p]{t_{1,x}} \isomto \prod_{x\in X} \Kxn[p]{t_{2,x}}$.  By Proposition~\ref{P:valuationsAndIsomorphisms}, $\phi$ restricts to an isomorphism $\AXn[p]{\t_{1}} \to \AXn[p]{\t_{2}}$ since $\upsilon_{x}( t_{1,x}) = \upsilon_{x}(t_{2,x}) =0$ for almost all $x\in X$. 

\item $(2) \implies (1)$: This follows from Lemma~\ref{L:familiesOfIsomorphisms} and Proposition~\ref{P:local conjugacy Kxt}.
\qedhere
\end{itemize}
\end{proof}

\begin{theorem}
\label{T:local global conjugacy}
For $\t_{i} \in \I_{X}$, $i=1,2$, let $(\AXn[p]{\t_{i}},G_{i})$ with $G_{i} \subseteq \G(\t_{i})$
be $p$-cyclic Galois ring extensions of $\A_{X}$. Then the following are equivalent:
\begin{enumerate}

\item $(\Kxn[p]{t_{1,x}},G_{1})$ and $(\Kxn[p]{t_{2,x}},G_{2})$ are conjugate Galois extensions via $(\phi_{x}, \tau_{x} )$ for all $x \in X$, and there is an isomorphism $\tau:G_{1}\isomto G_{2}$ making the diagram
\begin{equation}
\label{E:conjugate Gt1 Gt2 global}	
\begin{tikzcd}
 	G_{1}  \arrow[r, hook]  \arrow[d,"\tau", "\wr"'] & \G_{x}(\t_{1}) \arrow[d,"c_{\phi_{x}}","\wr"'] 
 	\\
 	G_{2}  \arrow[r, hook] & \G_{x}(\t_{2}) 
\end{tikzcd}
\end{equation} 
commutative where $c_{\phi_{x}}$ is conjugation by $\phi_{x}$.

\item $(\AXn[p]{\t_{1}}, G_{1})$ and $(\AXn[p]{\t_{2}},G_{2})$ are conjugate $p$-cyclic Galois extensions of $\A_{X}$.
\end{enumerate}

\end{theorem}

\begin{proof}
\leavevmode
\begin{itemize}

\item $(1) \implies (2)$: By Proposition~\ref{P:local conjugacy AXt}, $\phi:=(\phi_{x})$ induces an isomorphism $\AXn[p]{\t_{1}} \isomto \AXn[p]{\t_{2}}$. Noting that on the $x$-component it holds that $\tau_{x}= c_{\phi_{x}}\vert_{G_{1}} = \tau$, it follows that $(\AXn[p]{\t_{1}}, G_{1})$ and $(\AXn[p]{\t_{2}},G_{2})$ are conjugate via $(\phi,\tau)$. 

\item $(2) \implies (1)$: Reversing the reasoning, the claim follows.
\qedhere
\end{itemize}
\end{proof}

Note that the commutativity of~\eqref{E:conjugate Gt1 Gt2 global} is equivalent to saying that the pointwise conjugations $c_{\phi_{x}}\vert_{G_{1}}$ in fact do not depend on $x$.

% (end)
% subsection local global correspondence end

%%%%%%%%%%%%%%%%%%%%%%%%%%%%%%%%%%%%%%%%%%%%%%%%%%%%%%%%%%%%%%%%%%%%%%%%%%%%%%%%%%%
%%%%%%%%%%%%%%%%%%%%%%%%%%%%%%%%%%%%%%%%%%%%%%%%%%%%%%%%%%%%%%%%%%%%%%%%%%%%%%%%%%%
\subsection{Galois equivalent subgroups of automorphisms}
\label{subsec:galois equivalent subgroups of automorphisms}
% (fold)
%%%%%%%%%%%%%%%%%%%%%%%%%%%%%%%%%%%%%%%%%%%%%%%%%%%%%%%%%%%%%%%%%%%%%%%%%%%%%%%%%%%
%%%%%%%%%%%%%%%%%%%%%%%%%%%%%%%%%%%%%%%%%%%%%%%%%%%%%%%%%%%%%%%%%%%%%%%%%%%%%%%%%%%

In \S\ref{subsec:classification of cc_ p galois extensions} we have seen that in each isomorphism class of $p$-cyclic Galois extensions of $\A_{X}$ we always have an adelic representative $S = \AXn[p]{\t}$. 

We now focus on the interaction between the classification of adelic $p$-cyclic Galois extensions $(\AXn[p]{\t},G)$ up to conjugacy, where $G$ lies in the full automorphism group $\G(\t)$, and purely group-theoretical properties of such $G$.

First we need to determine which $p$-cyclic subgroups $G$ of the full automorphism group $\G(\t)$, yield $G$-Galois extensions $(\Axn[p]{\t}, G)$. This is done in Proposition~\ref{P:pointwise transitive equivalences} below.

The general problem of giving a group-theoretical characterization of when two $p$-cyclic Galois extensions $(\AXn[p]{\t},G_{1})$ and $(\AXn[p]{\t},G_{2})$ are conjugate is completely solved in Theorems~\ref{T:conjugacy subgroups equivalences} and~\ref{T:correspondence PT subgroups tuples Kummer} below. The latter theorem introduces finite tuples supported at the ramification points of $\t$. Using the local Kummer symbols,  Proposition~\ref{P:G-primitive elements valuation vector} gives an explicit relation between these tuples and the valuation vector $\upsilon(\AXn[p]{\t},G,\chi)$ introduced in~\eqref{E:(G,chi)-valuation vector}.

% We begin by determining which $p$-cyclic subgroups $G$ of the full automorphism group $\G(\t)$, yield $G$-Galois extensions $(\Axn[p]{\t}, G)$.

\begin{proposition}
\label{P:pointwise transitive equivalences}
Given an idele $\t \in \I_{X}$ and a $p$-cyclic subgroup $G$ of $\G(\t)$, the following are equivalent.
\begin{enumerate}
	
\item $\AXn[p]{\t}$	is a $G$-Galois ring extension of $\A_{X}$.

\item $\AXn[p]{\t}^{G} = \A_{X}$.

\item There is some $g \in G$ whose components $(g_{x})_{x} \in \prod_{x} \G_{x}(\t)$ under~\eqref{E:Gtcomponents} all have order $p$. Such elements are precisely the generators of $G$.

\item At each point $x \in X$, the projection of $G$ onto the $x$-component is a transitive subgroup of $\G_{x}(t)$, considering the latter as a permutation group on $\{1,2,\dots,p\}$ via~\eqref{E:Gtx structure}.

\end{enumerate}
\end{proposition}

\begin{proof}
\leavevmode
\begin{itemize}

\item $(1) \implies (2)$: Trivial.

\item $(2) \iff (3)$: Observe that $\AXn[p]{\t}^{\langle g \rangle} = \A_{X}$ if and only if $\Kxn[p]{t_{x}}^{\langle g_{x} \rangle} = K_{x}$ for every $x$. By~\eqref{E:Gtx structure}, this happens if and only if $g_{x}$ is a generator of $\cc_{p}$ if $x \in \Ram(\t)$ and a $p$-cycle in $\ss_{p}$ if $x \notin \Ram(\t)$. Clearly such elements $g$ are the generators of $G$ since it is assumed $p$-cyclic.

\item $(3) \iff (4)$: This is straightforward by well-known arguments from basic group theory.

\item $(4) \implies (1)$: Separability of $\A_{X} \to \AXn{\t}$ was shown in Proposition~\ref{P:Axt-separable}. Since we have already established the equivalence of (4) with (2), it remains only to show that the elements of $G$ are pairwise strongly distinct. Let $\tau_1,\tau_2\in G$ be distinct. 

Assume first that $x\notin \Ram(\t)$. Since $\Kxn[p]{t_{x}} \simeq \prod^{p} K_{x}$, to prove the strong distinctness, it suffices to consider the idempotents $e_{x,j}\in \AXn{\t}$ for $1 \leq j \leq p$, consisting of a $1$ in the $j$th factor and $0$ otherwise.

Since $G \simeq \cc_{p}$ has no nontrivial subgroups, the composition
\begin{equation}
\label{E:projectionGGxNoRam}
	G \hookrightarrow \G(\t) \to \G_{x}(\t) \simeq {\mathcal S}_p
\end{equation}
is either trivial or injective. By (2), it must be the latter, therefore its image is the subgroup generated by a $p$-cycle $\sigma$. Hence, there are distinct $n_1, n_2$ such that $(\tau_i)_{x} =\sigma^{n_i}$ via the above decomposition of $\Kxn[p]{t_{x}}$; that is
\[
	(\tau_i)_{x} (\alpha_{x}) = \sigma^{n_i} (a_1,\dots, a_p) = (a_{\sigma^{n_i}(1)} ,\dots, a_{\sigma^{n_i}(p)} ),
	\quad
	i = 1,2,
\]
for $\alpha_{x}=(a_1,\dots, a_p)\in \prod_{j=1}^p K_{x,j}$. It is clear that there exists  $\alpha_{x}$ with
		\[
	(\tau_1)_{x} (\alpha_{x}) e_{x,j} =  a_{\sigma^{n_1}(j)} \neq a_{\sigma^{n_2}(j)}  = (\tau_2)_{x} (\alpha_{x}) e_{x,j}.
	\]
If $x\in \Ram(\t)$, Proposition~\ref{P:pointwise transitive equivalences}(2) the composition
\begin{equation}
\label{E:proyectionGGxRam}
	G \hookrightarrow \G(\t) \to \G_{x}(\t) \simeq \cc_p
\end{equation}
is an isomorphism. Now, $\G_{x}(\t)$ is the Galois group of the Kummer extension $K_{x}\{t_{x}^{1/p}\}/K_{x}$, which acts freely and transitively on it.
\qedhere
\end{itemize}
\end{proof}

\begin{definition}[Pointwise transitive subgroup]
\label{D:pointwise transitive}
A $p$-cyclic subgroup $G$ of $\G(\t)$ satisfying the equivalent conditions of Proposition~\ref{P:pointwise transitive equivalences} will be called \emph{pointwise transitive}.
\end{definition}

\begin{definition}[Galois equivalence]
Let $\t \in \I_{X}$ and $G_{1}, G_{2} \subseteq \G(\t)$ two pointwise transitive subgroups. We say that $G_{1}$ and $G_{2}$ are Galois equivalent, $G_{1}\sim G_{2}$, if $(\AXn[p]{\t}, G_{1})$ and $(\AXn[p]{\t},G_{2})$ are conjugate $p$-cyclic Galois extensions of $\A_{X}$.
\end{definition}

Let us denote the image of $\G(\t)$ under the projection $\pi_{\Ram}$ onto the ramification components by $\G_{\Ram}(\t)$, that is,
\[
	       \G_{\Ram}(\t) 
	\eqdef \prod_{\mathmakebox[25pt]{x \in \Ram(t)}} \G_{x}(\t).
\]
As we have mentioned above, $\G_{\Ram}(\t)$ is commutative. In fact, for an odd prime $p$, it is easy to check that $\G_{\Ram}(\t) = \G(\t)^{\text{ab}}$, the abelianization of $\G(\t)$.

\begin{theorem}
\label{T:conjugacy subgroups equivalences}
Let $\t \in \I_{X}$ and $G_{1}, G_{2} \subseteq \G(\t)$ two pointwise transitive subgroups. Then the following are equivalent:
\begin{enumerate}

\item $G_{1}, G_{2}$ are Galois equivalent.

\item The projections of $G_{1}$ and $G_{2}$ onto the ramification components coincide, i.e. $\pi_{\Ram}(G_{1}) = \pi_{\Ram}(G_{2})$, where
\[
		\G(\t) = \prod_{x \in X} \G_{x}(\t)
		         \xrightarrow{\ \pi_{\Ram}\ }
		         \G_{\Ram}(\t) = \prod_{\mathmakebox[22pt]{x \in \Ram(t)}}\G_{x}(\t).
\]

\end{enumerate}
\end{theorem}

\begin{proof}
\leavevmode	

\begin{itemize}

\item $(1) \implies (2)$: Let $(\AXn[p]{\t}, G_{1})$ and $(\AXn[p]{\t},G_{2})$ be conjugate via $(\phi,\tau)$.
By Theorem~\ref{T:local global conjugacy}, $(\Kxn[p]{t_{x}},G_{1})$ and $(\Kxn[p]{t_{x}},G_{2})$ are conjugate Galois extensions via $(\phi_{x}, \tau_{x})$ for all $x \in X$ and the diagram~\eqref{E:conjugate Gt1 Gt2 global} commutes for $\t_{1}=\t_{2}:=\t$ and for all $x\in X$. For $x\in \Ram(\t)$, $\G_{x}(\t)\simeq \cc_{p}$ is commutative, so that $c_{\phi_{x}}=\id$. Hence, the result follows from the commutativity of the diagram
\begin{equation}
\label{E:conjugate Gt1 Gt2 global ram}
\begin{tikzcd}
 	G_{1}  \arrow[r, hook]  \arrow[d,"\tau", "\wr"'] & 
	\G(\t)  \arrow[r, "\pi_{\Ram}"]  \arrow[d,"(c_{\phi_{x}})_{x}","\wr"']  & 
	 \G_{\Ram}(\t) \arrow[d, equal]
 	\\
 	G_{2}  \arrow[r, hook] & 
	\G(\t)  \arrow[r, "\pi_{\Ram}"] & 
	 \G_{\Ram}(\t)
\end{tikzcd}
\end{equation} 

\item $(2) \implies (1)$: We will construct a pair $(\phi,\tau)$ conjugating $(\AXn[p]{\t}, G_{1})$ and $(\AXn[p]{\t}, G_{2})$. We begin by determining a suitable group isomorphism $\tau : G_{1} \isomto G_{2}$.

If $\Ram(\t) = \emptyset$, choose any isomorphism $\tau$. If $\Ram(\t) \neq \emptyset$, since $G_{1}$ is $p$-cyclic, the induced morphism $G_{1} \to \pi_{\Ram}(G_{1})$ is either trivial or an isomorphism. By Proposition~\ref{P:pointwise transitive equivalences}(3), it is the latter, and we define $\tau : G_{1} \simeq \pi_{\Ram}(G_{1}) = \pi_{\Ram}(G_{2}) \simeq G_{2}$, which yields the commutative square
\begin{equation}
\label{E:conjugate Gt1 Gt2 global ram exterior}
\begin{tikzcd}
 	G_{1}  \arrow[r, hook]  \arrow[d,"\tau", "\wr"'] & 
	\G_{\Ram}(\t) \arrow[d, equal]
 	\\
 	G_{2}  \arrow[r, hook] & 
	\G_{\Ram}(\t).
\end{tikzcd}
\end{equation} 

We will show that there is an automorphism $\phi \in \G(\t)$ such that conjugation by $\phi$ is $\tau$.

Choose a nontrivial character $\chi_{1} : G_{1} \to \mu_{p} \subseteq \k^{*}$ and a $(G_{1},\chi_{1})$-primitive element $\alpha_{1}$ of $(\AXn[p]{\t},G_{1})$. Consider the character of $G_{2}$ given by $\chi_{2} \eqdef \chi_{1} \circ \tau^{-1}$ and choose a $(G_{2},\chi_{2})$-primitive element $\alpha_{2}$ of $(\AXn[p]{\t},G_{2})$.

Recall that $\A_{X}$ is a $p$-Kummerian ring and that $\Pic(\A_{X}) = 0$ (Theorem~\ref{T:Pic adeles}). Thus by Theorem~\ref{T:primitive element theorem p-cyclic}, for $i=1,2$, we have the eigenspace decomposition of $\A_{X}$-modules
\[
	\AXn[p]{\t} = \bigoplus_{\chi \in \widehat{G_{i}}} S_{i}^{\chi},
\]
where each summand $S_{i}^{\chi}$ is a free $\A_{X}$-module of rank $1$. By definition, $\alpha_{i} \in S_{i}^{\chi_{i}}$ for $i=1,2$.

We will determine an idele $u \in \I_{X}$ such that sending $\alpha_{1}$ to $u \alpha_{2}$ yields the desired automorphism $\phi \in \G(\t)$, using Lemma~\ref{L:familiesOfIsomorphisms} to reduce the problem to points $x \in X$.

For any point $x \in X$ and any automorphism $\phi_{x} \in \G_{x}(\t)$, and $g \in G_{1}$, denoting its $x$-component under~\eqref{E:Gtcomponents} by $g_{x}$, we have
\begin{equation}
\label{E:conjugation phix}
	c_{\phi_{x}}(g_{x})(\phi_{x}(\alpha_{1,x})) = \chi_{1}(g) \phi_{x}(\alpha_{1,x})
\end{equation}
since
\[
	\begin{aligned}
	   \chi_{1}(g) \phi_{x}(\alpha_{1,x})
	&=  \phi_{x}((\chi_{1}(g) \alpha_{1})_{x})
	 = \phi_{x}((g(\alpha_{1}))_{x})
	 = \phi_{x}(g_{x}(\alpha_{1,x}))
	\\
	&= (\phi_{x} \circ g_{x} \circ \phi_{x}^{-1})(\phi_{x}(\alpha_{1,x}))
	 = c_{\phi_{x}}(g_{x})(\phi_{x}(\alpha_{1,x})).
	\end{aligned}
\]

\end{itemize}
\begin{itemize}

\item For $x \in \Ram(\t)$, taking into account the square~\eqref{E:conjugate Gt1 Gt2 global ram exterior}, which says that $c_{\phi_{x}}|_{G_{1}} = \tau$, and~\eqref{E:conjugation phix}, we obtain
\[
	 \begin{aligned}
	   \chi_{2}(\tau(g)) \phi_{x}(\alpha_{1,x})
	&= (\chi_{1} \circ \tau^{-1})(\tau(g)) \phi_{x}(\alpha_{1,x})
	 = \chi_{1}(g) \phi_{x}(\alpha_{1,x})
	\\
    &= c_{\phi_{x}}(g_{x})(\phi_{x}(\alpha_{1,x})) = \tau(g)_{x}(\phi_{x}(\alpha_{1,x}))
	 % = g_{2,x}(\phi_{x}(\alpha_{1,x})).
    \end{aligned}
\]
for any $g \in G_{1}$. Therefore, the image $\phi_{x}(\alpha_{1,x})$ is a $(G_{2},\chi_{2})$-primitive element of $\Kxn[p]{t_{x}}$ at each point $x \in \Ram(\t)$. Thus by Theorem~\ref{T:primitive element theorem p-cyclic}(5) we conclude that $\phi_{x}(\alpha_{1,x}) = u_{x} \alpha_{2,x}$ for some $u_{x} \in K_{x}^{*}$.

\item At an unramified point $x$, by Lemma~\ref{L:structureKxntx}, $\alpha_{i,x}$ corresponds to a vector $(a_{1,x}^{(i)}, a_{2,x}^{(i)}, \dots, a_{p,x}^{(i)}) \in \prod^{p} K_{x}$, considered as a trivial $K_{x}$-algebra (note however that this is not the base change of the eigenspace decomposition above).

Choose a generator $g_{1}$ of $G_{1}$ and let $g_{2} = \tau(g_{1})$, which generates $G_{2}$. By Proposition~\ref{P:pointwise transitive equivalences}, $g_{i,x}$ corresponds to a $p$-cycle $\sigma_{i} \in \ss_{p}$. Since  $\chi_{1}(g_{1}) = \chi_{2}(g_{2}) \defeq \zeta \in \mu_{p}$, the relation $g_{i,x}(\alpha_{i,x}) = \chi_{i}(g_{i}) \alpha_{i,x}$, in terms of the corresponding vectors, is
\[
	(a_{\sigma_{i}(1),x}^{(i)}, a_{\sigma_{i}(2),x}^{(i)}, \dots, a_{\sigma_{i}(p),x}^{(i)})
	= \zeta (a_{1,x}^{(i)}, a_{2,x}^{(i)}, \dots, a_{p,x}^{(i)})
\]
A straightforward argument shows that
\[
	  (a_{1,x}^{(i)}, a_{2,x}^{(i)}, \dots, a_{p,x}^{(i)})
	= a_{1,x}^{(i)} \zeta^{-1}
	  \bigl(\zeta^{\sigma_{i}^{0}(1)}, \zeta^{\sigma_{i}^{1}(1)}, \dots, \zeta^{\sigma_{i}^{p-1}(1)}\bigr).
\]
Since $\sigma_{i}$ is a $p$-cycle, the exponents of $\zeta$ in the above vector are a permutation of $1,2,\dots p$. Thus there is some $\sigma_{x} \in \ss_{p}$ which 
sends $(\zeta^{\sigma_{1}^{0}(1)}, \zeta^{\sigma_{1}^{1}(1)}, \dots, \zeta^{\sigma_{1}^{p-1}(1)})$ to $(\zeta^{\sigma_{2}^{0}(1)}, \zeta^{\sigma_{2}^{1}(1)}, \dots, \zeta^{\sigma_{2}^{p-1}(1)})$. This corresponds to an automorphism $\phi_{x} \in \G_{x}(\t)$ satisfying
\begin{align*}
	   \phi_{x}(\alpha_{1,x})
	&= \bigl(a_{1,x}^{(i)}, a_{2,x}^{(i)}, \dots, a_{p,x}^{(i)}\bigr)^{\sigma_{x}}
	\\
	&= \bigl(a_{1,x}^{(1)} \zeta^{-1}
	  \bigl(\zeta^{\sigma_{1}^{0}(1)}, \zeta^{\sigma_{1}^{1}(1)}, \dots, \zeta^{\sigma_{1}^{p-1}(1)}\bigr)\bigr)^{\sigma_{x}}
	\\
	&= a_{1,x}^{(1)} \zeta^{-1}
	  \bigl(\zeta^{\sigma_{2}^{0}(1)}, \zeta^{\sigma_{2}^{1}(1)}, \dots, \zeta^{\sigma_{2}^{p-1}(1)}\bigr)
	\\
	&= \frac{a_{1,x}^{(1)}}{a_{1,x}^{(2)}} a_{1,x}^{(2)} \zeta^{-1}
	  \bigl(\zeta^{\sigma_{2}^{0}(1)}, \zeta^{\sigma_{2}^{1}(1)}, \dots, \zeta^{\sigma_{2}^{p-1}(1)}\bigr)
	 = u_{x} \alpha_{2,x}
\end{align*}
for $u_{x} = a_{1,x}^{(1)}/a_{1,x}^{(2)}$.

\item We need to check that $u = (u_{x})_{x}$ is actually an idele. At each point we have $\alpha_{1,x}^{p} = u_{x}^{p} \alpha_{2,x}^{p}$, and recalling that $\alpha_{i}^{p} \in \I_{X}$ by Proposition~\ref{P:G-primitive equivalences}(2), we conclude that $\upsilon_{x}(u_{x}^{p}) = p \upsilon_{x}(u_{x}) = 0$ for almost all $x$, so that indeed $u \in \I_{X}$. Along with Proposition~\ref{P:valuationsAndIsomorphisms}, this shows that the family $(\phi_{x})_{x}$ defines an automorphism $\phi \in \G(\t)$ of $\AXn[p]{\t}$ which sends the $(G_{1},\chi_{1})$-primitive element $\alpha_{1}$ to the $(G_{2},\chi_{2})$-primitive element $u \alpha_{2}$.

\item Finally, to see that $(\AXn[p]{\t}, G_{1})$ and $(\AXn[p]{\t},G_{2})$ are conjugate via $(\phi,\tau)$, it suffices to check that $c_{\phi} = \tau$. We have,
\begin{align*}
	   (\chi_{1} \circ c_{\phi}^{-1})(\tau(g)) \phi(\alpha_{1})
	&= \phi( (\chi_{1} \circ c_{\phi}^{-1})(\tau(g)) \alpha_{1} )
	 = \phi( \chi_{1} (c_{\phi}^{-1}(\tau(g))) \alpha_{1} )
	\\
	&= \phi( c_{\phi}^{-1}(\tau(g)) \alpha_{1} )
	 = \phi( (\phi^{-1} \circ \tau(g) \circ \phi) \alpha_{1} )
	\\
	&= \tau(g) \phi(\alpha_{1})
	 = \chi_{2}(\tau(g)) \phi(\alpha_{1})
	\\
	&= (\chi_{1} \circ \tau^{-1}) (\tau(g)) \phi(\alpha_{1})
     = \chi_{1}(g) \phi(\alpha_{1})
\end{align*}
for all $g \in G_{1}$. Since $\phi(\alpha_{1}) = u \alpha_{2}$ is invertible, we conclude that
\[
	\chi_{1} \circ c_{\phi}^{-1} \circ \tau = \chi_{1}
\]
and, since $\chi_{1}$ is nontrivial, we conclude that indeed $c_{\phi}^{-1} \circ \tau = \id$.
\qedhere
\end{itemize}
\end{proof}

Theorem~\ref{T:conjugacy subgroups equivalences} leads us to focus our attention on the group $\G_{\Ram}(\t)$. Determining its structure leads naturally to another set of invariants, in this case finite tuples over ramified points. As we will see, these are related to the valuation vectors introduced in \S\ref{subsec:classification of cc_ p galois extensions} via some computations with the local Kummer symbols.

When $x \in \Ram(\t)$, Lemma~\ref{L:structureKxntx} shows that the $K_{x}$-algebra $\Kxn[p]{t_{x}}$ is a field isomorphic to the $p$-cyclic Kummer extension $E_{x} = K_{x}(t_{x}^{1/p})$.	Since $\G_{x}(\t) = \Gal(E_{x}/K_{x})$, by Lemma~\ref{L:Puiseux formal}, the Kummer pairing for $E_{x}$ induces a perfect pairing
\begin{equation}
\label{E:Kummer pairing Kx}
	\pairing{g,\lambda}_{x} = \frac{g(\lambda^{1/p})}{\lambda^{1/p}} :
	\G_{x}(\t) \times \faktor{K_{x}^{*}}{K_{x}^{*p}} \to \mu_{p},
\end{equation}
where $\lambda^{1/p}$ denotes any $p$th root of $\lambda$ in $E_{x}$.

\begin{lemma}
\label{L:GRam characterization}
There is a canonical group isomorphism
\begin{equation}
\label{E:GRam characterization}
	       \G_{\Ram}(\t) 
	% =      \prod_{\mathmakebox[25pt]{x \in \Ram(t)}} \G_{x}(\t)
	\simeq \bigoplus_{\mathmakebox[25pt]{x\in\Ram(\t)}}  \widehat{ \faktor{\Z}{(p)}}.
	% \simeq \prod_{\mathmakebox[25pt]{x\in\Ram(\t)}}  \faktor{\Z}{(p)}.
\end{equation}
\end{lemma}

\begin{proof}
Any $g\in \G_{x}(\t)$ yields, via the Kummer pairing~\eqref{E:Kummer pairing Kx}, a character $\chi_{g}(\lambda) = \pairing{g,\lambda}_{x}$ of $K_{x}^{*}/K_{x}^{*p}$. Since $A_{x}^{*} \subseteq K_{x}^{*p}$, this character depends only on $\upsilon_{x}(\lambda) \bmod p$, and the map $g \mapsto \chi_{g}$ induces an isomorphism
\[
	\G_{x}(\t) \simeq  
	\widehat{ \faktor{K_{x}^{*}}{K_{x}^{*p}} } \overset{\upsilon_{x}}\simeq
	\widehat{ \faktor{\Z}{(p)} }.
\] 
Considering the product over all the ramification points we obtain~\eqref{E:GRam characterization}. 
\end{proof}

\begin{theorem}
\label{T:correspondence PT subgroups tuples Kummer}
Given an idele $\t \in \I_{X}$, there is a canonical set correspondence
\begin{equation}
\label{E:correspondence classes subgroups}
	\left\{
	\begin{minipage}[c]{0.36\textwidth}
	\raggedright
	Galois equivalence classes of pointwise transitive subgroups of $\G(\t)$
	\end{minipage}
	\right\}
	\correspondence
	\begin{minipage}[c]{0.25\textwidth}
	\raggedright
	\[\prod_{x \in \Ram(\t)} \!\! \Zpsf \]
	% s.t. $\sum v_i \equiv 0 \bmod p$
	\end{minipage}	
	\!\!\!\Bigm/\sim,
\end{equation}
where the equivalence relation on tuples is given by
\begin{equation}
\label{E:equivalence relation tuples}
		(v_{x}^{1})_{x \in \Ram(\t)} \sim (v_{x}^{2})_{x \in \Ram(\t)}
\end{equation}
if there exists some $b \in \Zps$ such that $b v_{x}^{1} \equiv v_{x}^{2}$ for $x \in \Ram(\t)$.
\end{theorem}

\begin{proof}
Fix any isomorphism between $\Zp$ and its dual. Then the group isomorphism~\eqref{E:GRam characterization} yields a bijection between sets
\begin{equation}
\label{E:Gram product Cp}
	\G_{\Ram}(\t) \simeq \prod_{x \in \Ram(\t)} \faktor{\Z}{(p)}.
\end{equation}
Now, let $G$ be a pointwise transitive subgroup of $\G(\t)$. By Theorem~\ref{T:conjugacy subgroups equivalences} its Galois equivalence class is determined by its projection onto $\G_{\Ram}(\t)$. Hence, if $g\in G$ is a generator, then \eqref{E:Gram product Cp} associates a tuple to $\pi_{\Ram}(g)$ whose components, by Proposition~\ref{P:pointwise transitive equivalences}, are generators of $\Zp$, i.e. lie in $\Zps$. Clearly choosing another generator of $G$ yields an equivalent tuple. 

Conversely, let us see that every tuple $(v_{x})_{x \in \Ram(\t)} \in \prod_{x \in \Ram(\t)} \Zps$ corresponds to a  pointwise transitive subgroup of $\G(\t)$. Via~\eqref{E:Gram product Cp}, this is equivalent to showing that for every $(g_x)_{x \in \Ram(\t)} \in \G_{\Ram}(\t)$, where each $g_x$ is a generator of $\G_{x}(\t)$, there is a pointwise transitive subgroup $G$ of $\G(\t)$ and a generator $g\in G$ such that $\pi_{\Ram}(g)=(g_x)$. Indeed, consider the element $g\in \G(\t)$ defined by $g_x$ at $x\in \Ram(\t)$ and by any $p$-cycle at $x\notin \Ram(\t)$ (via \eqref{E:Gtx structure}). It is clear that $G := \langle g\rangle\subset \G(\t)$ is pointwise transitive by Proposition~\ref{P:pointwise transitive equivalences} and corresponds to the given tuple.

Finally, note that choosing another isomorphism between $\Zp$ and its dual changes a given tuple to an equivalent one in~\eqref{E:Gram product Cp}, so that the correspondence is canonical.
\end{proof}

Combining~\eqref{E:correspondence conjugacy p-cyclic mod aut},~\eqref{E:correspondence conjugacy p-cyclic AX} and~\eqref{E:correspondence classes subgroups} we obtain the following.

\begin{corollary}[Stratification by ramification]
\label{C:stratification conjugacy p-cyclic ramification}
\begin{equation}
\label{E:stratification conjugacy p-cyclic ramification}
	\begin{aligned}
	\faktor{\H(\A_{X},\cc_{p})}{\Aut(\cc_{p})}
	&\correspondence
	\left\{
	\begin{minipage}[c]{0.428\textwidth}
	\raggedright
	Conjugacy classes of $p$-cyclic Galois extensions $(S,G)$ of $\A_{X}$
	\end{minipage}
	\right\}
	\\
	&\correspondence
	\coprod_{\mathfrak{R} \subseteq X} \Bigl(\, \prod_{x \in \mathfrak{R}} \Zps \Bigr) \!\! \Bigm/ \sim
	\end{aligned}
\end{equation}
where $\mathfrak{R}$ ranges over finite subsets of closed points of $X$ and $\sim$ denotes the equivalence relation~\eqref{E:equivalence relation tuples}.
\end{corollary}

The previous considerations lead to the following explicit construction of tuples from pointwise transitive subgroups. Fixing a primitive $p$th root of unity $\zeta \in \mu_{p}$ defines an isomorphism between $\Zp$ and its dual given by mapping $\chi \in (\Zp)^{\wedge}$ to the residue class $c$ such that $\chi(1) = \zeta^{c}$. Denoting the inverse of $c \mapsto \zeta^{c}$ by the discrete logarithm $\log_{\zeta} : \mu_{p} \to \Zp$, we have $c = \log_{\zeta} \chi(1)$.

An element $g = (g_{x}) \in \G_{\Ram}(\t)$ gives, for each $x \in \Ram(\t)$, a character $\lambda \mapsto \pairing{g_{x},\lambda}_{x}$ of $K_{x}^{*}/K_{x}^{*p}$ via the Kummer pairing. By~\eqref{E:Kummer Kx/Kxp=Z/p} this corresponds to the character of $\Zp$ given by $a \mapsto \pairing{g_{x},\lambda_{a}}_{x}$, where $\lambda_{a} \in K_{x}^{*}$ is any element with $\upsilon_{x}(\lambda_{a}) = a$. In fact this character is determined by its action on the class of $1$, which corresponds to a choice of uniformizer $z_{x}$ at $x$, i.e.
\begin{equation}
\label{E:character Zp uniformizer}
	1 \in \Zp \longmapsto \pairing{g_{x},z_{x}}_{x} \in \mu_{p}.
\end{equation}
Now, via our choice of primitive root, the character is mapped to the residue class $c_{x} \in \Zp$ such that $\pairing{g_{x},z_{x}}_{x} = \zeta^{c_{x}}$. Thus, the isomorphism~\eqref{E:Gram product Cp} is explicitly given by
\begin{equation}
\label{E:r-tuple explicit}
	g \in \G_{\Ram}(\t)
	\longmapsto
	\bigl(\log_{\zeta}\pairing{g_{x},z_{x}}_{x}\bigr)_{x \in \Ram(\t)} \in \prod_{x \in \Ram(\t)} \Zp.
\end{equation}
The dependence of this isomorphism on the choice of primitive root $\zeta$ is governed by the action of $\Aut(\Zp)$.

If $G$ is a pointwise transitive subgroup of $\G(\t)$, we have constructed two associated sets of data, namely, the tuples defined in Theorem~\ref{T:correspondence PT subgroups tuples Kummer}, explicitly given in terms of the Kummer symbol by~\eqref{E:r-tuple explicit}, and a valuation vector $\upsilon(\AXn[p]{\t},G,\chi)$ as in~\eqref{E:(G,chi)-valuation vector}. These sets are related as follows.

% , where the corresponding entries are essentially inverses modulo $p$.

\begin{proposition}
\label{P:G-primitive elements valuation vector}
Given an idele $\t \in \I_{X}$, a pointwise transitive subgroup $G$ of $\G(\t)$, a nontrivial character $\chi \in \widehat{G}$, a generator $g$ of $G$, a primitive $p$th root of unity $\zeta$, and a $(G, \chi)$-primitive element $\alpha$ of $\AXn[p]{\t}$ (which exists by Theorem~\ref{T:primitive element theorem p-cyclic}), then in terms of the Kummer pairing~\eqref{E:Kummer pairing Kx} we have the explicit formula
\begin{equation}
\label{E:valuation vector ap computation}
	       \upsilon_{x}(\alpha^{p})
	\equiv \begin{dcases*}
		   0 \bmod p
		   & if and only if $x \notin \Ram(\t)$,
		   \\
		   \frac{\log_{\zeta} \chi(g)}{\log_{\zeta}\pairing{g_{x},z_{x}}_{x}} \bmod p
		   & if $x \in \Ram(\t)$,
		   \end{dcases*}	
\end{equation}
independently of the choice of $g$ and $\zeta$. Here $(g_{x})$ denotes the family of automorphisms corresponding to $g$ via Lemma~\ref{L:familiesOfIsomorphisms}. %In particular, $\Ram(\t) = \Ram(\alpha^{p})$.
\end{proposition}

\begin{proof}
As we saw in the discussion following~\eqref{E:valuation vector (G,chi)-Kummer}, if $\alpha$ is a $(G,\chi)$-primitive element, the valuation vector $\upsilon(\AXn[p]{\t},G,\chi)$ is given by $\upsilon(\alpha^{p})$ and, by~\eqref{E:ram ap = ram t}, it is only nonzero modulo $p$ at the ramified points $x \in \Ram(\t)$. We can explicitly determine the valuations at these points as follows.

For $x \in \Ram(\t)$, we have $\alpha_{x} \in E_{x}$, the Kummer extension considered in Lemma~\ref{L:Puiseux formal}, and $\alpha_{x}^{p} \in K_{x}^{*}$. Thus the Kummer pairing~\eqref{E:Kummer pairing Kx} is given by
\begin{equation}
\label{E:Kummer pairing G-primitive}
		  \pairing{g_{x},z_{x}}_{x}^{\upsilon_{x}(\alpha_{x}^{p})}
	= \pairing{g_{x},\alpha_{x}^{p}}_{x}
	= \frac{g_{x}(\alpha_{x})}{\alpha_{x}}
	= \chi(g).
\end{equation}
Recall that by Proposition~\ref{P:pointwise transitive equivalences}(3), $g_{x}$ is a generator of $\G_{x}(\t)$, hence $\pairing{g_{x},z_{x}}_{x} \neq 1$, and~\eqref{E:valuation vector ap computation} follows. Clearly the quotient is independent of the choices of $g$ and $\zeta$.
\end{proof}

% (end)
% subsection galois equivalent subgroups of automorphisms

% (end)
% section Kummer extensions of AX

%%%%%%%%%%%%%%%%%%%%%%%%%%%%%%%%%%%%%%%%%%%%%%%%%%%%%%%%%%%%%%%%%%%%%%%%%%%%%%%%%%%
%%%%%%%%%%%%%%%%%%%%%%%%%%%%%%%%%%%%%%%%%%%%%%%%%%%%%%%%%%%%%%%%%%%%%%%%%%%%%%%%%%%
%%%%%%%%%%%%%%%%%%%%%%%%%%%%%%%%%%%%%%%%%%%%%%%%%%%%%%%%%%%%%%%%%%%%%%%%%%%%%%%%%%%
\section{Future Work}
\label{sec:future work}
% (fold)
%%%%%%%%%%%%%%%%%%%%%%%%%%%%%%%%%%%%%%%%%%%%%%%%%%%%%%%%%%%%%%%%%%%%%%%%%%%%%%%%%%%
%%%%%%%%%%%%%%%%%%%%%%%%%%%%%%%%%%%%%%%%%%%%%%%%%%%%%%%%%%%%%%%%%%%%%%%%%%%%%%%%%%%
%%%%%%%%%%%%%%%%%%%%%%%%%%%%%%%%%%%%%%%%%%%%%%%%%%%%%%%%%%%%%%%%%%%%%%%%%%%%%%%%%%%

We now give a sketch of how the theory developed in this paper will be applied to the study of abelian extensions of the function field $\Sigma$ of the curve $X$ or, equivalently, to its abelian Galois covers. This is of course a classical topic which has been approached from a variety of points of view and techniques. One only has to think of the Galois Theory of finite field extensions or the Class Field Theory of local and global fields. It continues to be the subject of current research as part of broader projects such as the equivalence problem for fields (see e.g. \cite{Sutherland} for a recent survey), the Inverse Galois Problem and the study of covers of a given curve.

Here, using the general Galois Theory of commutative ring extensions, we have developed in detail the Kummer theory of the adele ring $\A_{X}$ (introduced in ~\S\ref{subsec:adeles}) of the function field of an algebraic curve. In particular we have completely determined the structure of the Harrison group $\H(\A_{X},\cc_{p})$ for the cyclic group $\cc_{p}$ of prime order $p$. Along the way, we have exhibited various instances of local-global relations (for instance, Theorem~\ref{T:galois local global}).

The functoriality of the Harrison group (\cite[Proposition 3.10]{Greither}) allows us to relate the Kummer theory of the function field $\Sigma$ with that of its adele ring $\A_{X}$. To be precise, the map induced by tensoring with $\A_{X}$ induces a group homomorphism $\H(\Sigma,G) \to \H(\A_{X}, G)$, where $G$ is a finite abelian group (we chose $G = \cc_{p}$ for simplicity but the method is generalizable). Studying the kernel and image of this map will allow us to describe the $p$-cyclic extensions of $\Sigma$ in terms of adeles. This idea may be said to be follow the work of Artin-Whaples~\cite{ArtinW},  Iwasawa~\cite{Iwa}, Tate~\cite{Tate-endo} and Weil~\cite{Weil-corps} among others. In this regard, it is also worth noting that the geometric adele ring has been used to study reciprocity laws in~\cite{AndersonPablos, collectanea} and~\cite{MP,MP-mediterranean}.

To illustrate this, consider the following situation. Let $p$ be a prime different from $\chr\k$. Let $f \in \Sigma^{*}$ be a rational function that is not a $p$th power in $\Sigma$. Then the polynomial $T^{p} - f$ is irreducible over $\Sigma$ and, since $\k$ is assumed to be algebraically closed, $\Omega = \Sigma[T]/(T^{p} - f)$ is a $p$-cyclic Kummer extension of $\Sigma$, namely $\Gal(Y/X) \simeq \cc_{p}$, where $Y$ is the Zariski-Riemann variety of $\Omega$.

In this context, we have the following result, which is the analog of~\cite[Lemma 8]{Iwa}, with linear local compactness replaced by the topology of commensurability as defined in~\cite{collectanea}:
\begin{equation}
\label{E:cassels}
	\A_{X}\otimes_{\Sigma} \Omega  \simeq   \A_{Y}	
\end{equation}
as linear topological $\A_{X}$-algebras.
% Furthermore, $\A_{Y}^{+}$ is the integral closure of $\A_{X}^{+}$ in $\A_{X}\otimes_{\Sigma} \Omega$.
Furthermore, this isomorphism is compatible with the action of $\Gal(Y/X)$, and therefore
\begin{equation}
\label{E:AY Gal(Y/X) invariant is AX}
	\A_{Y}^{\Gal(Y/X)} = \A_{X}.
\end{equation}
%%

% \begin{theorem}
% \label{T:cassels}
% Let $\Omega$ be a finite separable field extension of $\Sigma$ and let $\A_{Y}$ denote the ring of adeles of $\Omega$. Then
% %%
% \begin{equation}
% \label{E:cassels}
% 	\A_{X}\otimes_{\Sigma} \Omega  \simeq   \A_{Y}
% \end{equation}
% %%
% as linear topological $\A_{X}$-algebras. Furthermore, $\A_{Y}^{+}$ is the integral closure of $\A_{X}^{+}$ in $\A_{X}\otimes_{\Sigma} \Omega$.
% \end{theorem}

% In the case of a Galois cover $Y \to X$, that is to say, with $\Omega/\Sigma$ a finite Galois field extension, by considering~\eqref{E:tensorCompletion} at each fiber, we can easily verify that the isomorphism $\A_{Y} \simeq \A_{X} \otimes_{\Sigma} \Omega$ from~\eqref{E:cassels} is compatible with the action of $\Gal(Y/X)$, and therefore
% %%
% \begin{equation}
% \label{E:AY Gal(Y/X) invariant is AX}
% 	\A_{Y}^{\Gal(Y/X)} = \A_{X}.
% \end{equation}
% %%

By the Kummer theory of fields, every $p$-cyclic field extension of $\Sigma$ is of this form (\cite[Ch. III, \S2]{CaFr}). In addition, the action of $\cc_{p}$ is given by a character, namely $g(T) = \chi(g)T$ for some $\chi : \cc_{p} \to \mu_{p} \subseteq \k^{*}$. This is what we have termed a $(\cc_{p},\chi)$-Kummer extension (\S\ref{subsec:structure of p cyclic extensions}). In this way, we obtain an element of $\H(\Sigma,\cc_{p})$.

We may also consider the adelic algebra $\AXn[p]{\f} = \A_{X}[T]/(T^{p} - \f)$ where $\f = (f_{x})_{x}$ is the idele such that $f_{x} \in K_{x}$ is the germ of $f$ at $x$. Considering the action of $\cc_{p}$ via the same character $\chi$, we obtain a $(\cc_{p},\chi)$-Kummer extension of $\A_{X}$ as described in \S\ref{subsec:classification of cc_ p galois extensions}, and hence an element of $\H(\A_{X},\cc_{p})$.

The following well-known example shows in more detail the role of ramification and how it fits in with our construction~\eqref{E:(G,chi)-valuation vector} of valuation vectors.

\begin{example}[Superelliptic curves]
\label{EX:ramified cover P1}
Consider the following classical construction of a $p$-cyclic cover of the projective line. We take $\k= \C$, $X= \P_1$, a prime number $p$, and $Y$ the normalization of 
\[
	y^{p}  =  f(x) \eqdef ( x- x_1)^{v_1} \cdots (x-x_r)^{v_r},
\]
where $\pi:Y\to X$ maps $(x,y)$ to $x$, the $x_{i}$ are distinct, and the exponents $v_{i}$ satisfy $0 < v_{i} < p$ and $\sum_{i} v_{i} \equiv 0 \bmod p$. Observe that the Riemann-Hurwitz formula implies that $r \geq 1$ in the above expression for $f(x)$.

Consider a character $\chi : \cc_{p} \to  \mu_{p}$ defined by a choice of primitive root of unity $\zeta \eqdef \chi(1) \in \mu_{p}$ and the action of $\cc_{p}$ on $Y$ by multiplication by $\zeta$ on the $y$-component of a point $(x,y)\in Y$. If we assume that  $\operatorname{gcd}(v_1,\ldots, v_r)=1$, then $Y$ is irreducible and $\oo_X\to \pi_*\oo_Y$ at the generic point is
\[
	\Sigma \eqdef \C(x) \hookrightarrow  \Omega \eqdef \Sigma[y]/(y^p-f(x)).
\]
Then, $\pi: Y\to X$ is ramified at $\Ram(\pi) \eqdef \{x_1, \ldots, x_r\}$, and the ramification degree at $x_i$ is ${p}/{(p,v_i)}=p$. 

The corresponding adelic algebra $\A_{Y} = \AXn[p]{\f}$ is a  $(\cc_{p},\chi)$-Kummer adelic extension of $\A_{X}$ as in the previous discussion, with characteristic polynomial of $y$ equal to $C_y(T)=T^p - f(x)$. Finally, its associated valuation vector is
\[
	  \upsilon(\A_{Y},\cc_{p},\chi)
	= (\upsilon_{x}(f(x)))_{x \in X}
	= \begin{cases}
	  v_{i} \bmod p, & x = x_{i} \in \Ram(\pi),
		    \\
	  0	    \bmod p, & x \notin \Ram(\pi).
	  \end{cases}
\]
\end{example}

By~\eqref{E:cassels} and~\eqref{E:AY Gal(Y/X) invariant is AX}, $\A_{Y}$ is the ring extension obtained by tensoring the field extension $\Omega$ with $\A_{X}$, namely, it is its image under the functorial homomorphism $\H(\Sigma, \cc_{p}) \to \H(\A_{X}, \cc_{p})$. In a forthcoming paper~\cite{adeles03} we will study the kernel and image of this map and its relation with ramification. This will achieve our goal of classifying covers of $X$ by looking among extensions of $\A_{X}$.

% (end)
% end section future work

%%%%%%%%%%%%%%%%%%%%%%%%%%%%%%%%%%%%%%%%%%%%%%%%%%%%%%%%%%%%%%%%%%%%%%%%%%%%%%%%%%%
%%%%%%%%%%%%%%%%%%%%%%%%%%%%%%%%%%%%%%%%%%%%%%%%%%%%%%%%%%%%%%%%%%%%%%%%%%%%%%%%%%%
%%%%%%%%%%%%%%%%%%%%%%%%%%%%%%%%%%%%%%%%%%%%%%%%%%%%%%%%%%%%%%%%%%%%%%%%%%%%%%%%%%%
%%%%%%%%%%%%%%%%%%%%%%%%%%%%%%%%%%%%%%%%%%%%%%%%%%%%%%%%%%%%%%%%%%%%%%%%%%%%%%%%%%%
%%%%%%%%%%%%%%%%%%%%%%%%%%%%%%%%%%%%%%%%%%%%%%%%%%%%%%%%%%%%%%%%%%%%%%%%%%%%%%%%%%%
%%%%%%%%%%%%%%%%%%%%%%%%%%%%%%%%%%%%%%%%%%%%%%%%%%%%%%%%%%%%%%%%%%%%%%%%%%%%%%%%%%%

\end{document}